\documentclass[a4paper,12pt]{article}
\usepackage{amsmath, amsfonts, amssymb,latexsym}
\usepackage{eurosym}
\usepackage[latin1]{inputenc}
\usepackage{graphicx}

\newtheorem{coro}{Corollary}[section]
\newtheorem{prop}{Proposition}[section]
\newtheorem{thm}{Theorem}[section]
\newtheorem{lem}{Lemma}[section]
\newtheorem{defi}{Definition}[section]

\numberwithin{equation}{section}

\begin{document}

\title{\textbf{\ }On\textbf{\ }the well posedness of a mathematical model
for a singular nonlinear fractional pseudo-hyperbolic system with nonlocal
boundary conditions and frictional dampings}
\author{Said Mesloub$^{1}$ Hassan Eltayeb Gadian$^{1}$, and Lotfi Kasmi$^{2}$
\\
$^1$Mathematics Department, College of Science,\\
King Saud University, PO Box2455, Riyadh, Saudi Arabia.\\
mesloub@ksu.edu.sa\\
[2pt] $^{2}$ Applied Mathematics Lab University\ Kasdi Merbah \ Ouargla,
Algeria. \\
kasmi.lotfi39@gmail.com}
\maketitle

\begin{abstract}
This paper is devoted to the study of the well-posedness of a singular
nonlinear fractional pseudo-hyperbolic system. The fractional derivative is
described in Caputo sense. The equations are supplemented by classical and
nonlocal boundary conditions. Upon some a priori estimates and density
arguments, we establish the existence and uniqueness of the strongly
generalized solution for the associated linear fractional system in some
Sobolev fractional spaces. On the basis of the obtained results for the
linear fractional system, we apply an iterative process in order to
establish the well-posedness of the nonlinear fractional system. This
mathematical model of pseudo-hyperbolic systems arises mainly in the theory
of longitudinal and lateral vibrations of elastic bars (beams), and in some
special case it is propounded in unsteady helical flows between two infinite
coaxial circular cylinders for some specific boundary conditions. \newline

\end{abstract}

\section{Introduction}

In the bounded domain $Q=\Omega \times (0,T)=\{(x,t):0<x<b,0<t<T\},$we are
concerned with the well posedness of a nonlinear fractional system with
frictional damping. More precisely, the model problem we have in mind is
presented in the form
\begin{equation}
\begin{cases}
& ^{C}\partial _{0t}^{\beta }u-\frac{1}{x}\left( xu_{x}\right) _{x}-\frac{%
\partial }{\partial t}\frac{1}{x}\left( xu_{x}\right)
_{x}+z_{1}v+u_{t}=f\left( x,t,u,v,u_{x},v_{x}\right) , \\
& ^{C}\partial _{0t}^{\gamma }v-\frac{1}{x}\left( xv_{x}\right) _{x}-\frac{%
\partial }{\partial t}\frac{1}{x}\left( xv_{x}\right)
_{x}+z_{2}u+v_{t}=g\left( x,t,u,v,u_{x},v_{x}\right) , \\
& u(x,0)=\varphi _{1}(x),~~~u_{t}(x,0)=\varphi _{2}(x), \\
& v(x,0)=\psi _{1}(x),~~~v_{t}(x,0)=\psi _{2}(x), \\
& u_{x}(b,t)=0,~~~v_{x}(b,t)=0,~~~\int\limits_{0}^{b}xudx=0,~~~\int%
\limits_{0}^{b}xvdx=0.%
\end{cases}
\label{e1.1}
\end{equation}%
The functions $f,$ $g$ are $L^{2}(0,T;L_{\rho }^{2}(\Omega )$ given
Lipschitzian functions, that is there exist two positive constants $\delta
_{1}$, $\delta _{2}$ such that%
\begin{eqnarray*}
&&\left\vert
f(x,t,u_{1},v_{1},w_{1},d_{1})-f(x,t,u_{2},v_{2},w_{2},d_{2})\right\vert \\
&\leq &\delta _{1}(\left\vert u_{1}-u_{2}\right\vert +\left\vert
v_{1}-v_{2}\right\vert +\left\vert w_{1}-w_{2}\right\vert +\left\vert
d_{1}-d_{2}\right\vert ),
\end{eqnarray*}%
\begin{eqnarray*}
&&\left\vert
g(x,t,u_{1},v_{1},w_{1},d_{1})-f(x,t,u_{2},v_{2},w_{2},d_{2})\right\vert \\
&\leq &\delta _{2}(\left\vert u_{1}-u_{2}\right\vert +\left\vert
v_{1}-v_{2}\right\vert +\left\vert w_{1}-w_{2}\right\vert +\left\vert
d_{1}-d_{2}\right\vert ),
\end{eqnarray*}%
for all $(x,t)\in Q$. The functions $\varphi _{1},$ $\psi _{1}$, $\varphi
_{2}$ and $\psi _{2}$ are in $H_{\rho }^{1}(\Omega ),$ and $z_{1},z_{2}$ are
positive constants. The operator $^{C}\partial _{0t}^{\beta }$ denotes the
left Caputo fractional derivative, defined in the second section, where $%
1<\beta ,\gamma <2$.

During the last twenty years, and up to this moment, the fractional order
differential equations are the essential tool in the modeling of several
phenomena in biology [$1,2,3,4,5,6,7,8,9,10,11,17$], in controlling chaotic
dynamical systems [12, 13, 14, 15, 16, 18, 41, 42 ], in heat transfer and
diffusion [19, 20, 21, 22, 23, 24, 25, 32, 36], in finance [26, 27, 28, 29,
30, 31 ], in thermoelasticity [33, 34, 35, 37, 38, 39, 40], in the
description of viscoelasticity [43, 44, 45], and some other different fields
such as mechanics, engineering, and seismology. The present mathematical
model of fractional pseudo-hyperbolic equations arises mainly in the theory
of longitudinal and lateral vibrations of elastic bars (beams). For the non
fractional case of hyperbolic and pseudo-hyperbolic see for example [46, 47,
48, 49, 50]. In [51], the authors studied a model which is propounded in the
investigation of the unsteady helical flows of a generalized Oldroyd-B fluid
with fractional calculus between two infinite coaxial circular cylinders
with initial conditions and Diruchlet boundary conditions%
\begin{equation*}
\begin{cases}
\lambda _{1}^{\alpha C}\partial _{0t}^{\alpha }u+u_{t}-\nu \frac{1}{x}\left(
xu_{x}\right) _{x}-\lambda _{2}^{\eta C}\partial _{0t}^{\eta }\frac{1}{x}%
\left( xu_{x}\right) _{x}=0 \\
\lambda _{1}^{\alpha C}\partial _{0t}^{\gamma }v+v_{t}-\frac{1}{x}\left(
xv_{x}\right) _{x}-\lambda _{2}^{\eta C}\partial _{0t}^{\eta }\frac{1}{x}%
\left( xv_{x}\right) _{xt}=0.%
\end{cases}%
\end{equation*}%
This model can be considered as a particular case of our model (\ref{e1.1})
with $z_{1}=0,$ $z_{2}=0,$ $\lambda _{1}^{\alpha }=\lambda _{2}^{\eta }=1,$ $%
\eta =1,f=g=0,$ and the Newmann and integral conditions were replaced by
Dirichlet conditions.

The paper is organized as follows: In Section 2, we introduce the needed
function spaces, and state some important inequalities, and fractional
calculus relations that will be used in the rest of the sequel. In Section
3, we reformulate the fractional linear system associated to the nonlinear
problem (\ref{e1.1}) in its operator form. Then in Section 4, we prove the
uniqueness of the solution of the fractional linear system, and we present
the consequences of the obtained energy estimate (\ref{e4.1}) of the
solution. In Section 5, we show the solvability of the associated linear
problem. Finally, in Section 6, on the basis of the results obtained in
Sections 4 and 5, and on the use of an iterative process, we prove the
existence and uniqueness of the solution of the fractional nonlinear system (%
\ref{e1.1}).

\section{Preliminaries and functions spaces}

\subsection{Functional spaces}

Let $L^{2}(0,T;L_{\rho }^{2}(\Omega ))$ be the space consisting of all
measurable functions $\mathcal{Q}:[0,T]\rightarrow L_{\rho }^{2}(\Omega )$
with scalar product
\begin{equation}
(\mathcal{Q},\mathcal{Q}^{\ast })_{L^{2}(0,T;L_{\rho }^{2}(\Omega
))}=\int\limits_{0}^{T}(\mathcal{Q},\mathcal{Q}^{\ast })_{L_{\rho
}^{2}(\Omega )}dt,  \label{e2.1}
\end{equation}%
and with the associated finite norm
\begin{equation}
\Vert \mathcal{Q}\Vert _{L^{2}(0,T;L_{\rho }^{2}(\Omega
))}^{2}=\int\limits_{0}^{T}\Vert \mathcal{Q}\Vert _{L_{\rho }^{2}(\Omega
)}^{2}dt,  \label{e2.2}
\end{equation}%
and we denote by $L^{2}(0,T;H_{\ \rho }^{1}(\Omega ))$ the space of
functions which are square integrable in the Bochner sense, with the inner
product%
\begin{equation}
(\mathcal{Q},\mathcal{Q}^{\ast })_{L^{2}(0,T;H_{\ \rho }^{1}(\Omega
))}=\int\limits_{0}^{T}(\mathcal{Q}(.,t),\mathcal{Q}^{\ast }\mathcal{(}%
.,t)_{H_{\ \rho }^{1}(\Omega )}dt,  \label{e2.3}
\end{equation}%
and the associated norm is
\begin{equation}
\Vert \mathcal{Q}\Vert _{L^{2}(0,T;H_{\ \rho }^{1}(\Omega
))}^{2}=\int\limits_{0}^{T}\Vert \mathcal{Q}(.,t)\Vert _{L_{\rho
}^{2}(\Omega ))}^{2}dt+\int\limits_{0}^{T}\Vert \mathcal{Q}_{x}(.,t)\Vert
_{L_{\rho }^{2}(\Omega ))}^{2}dt.  \label{e2.4}
\end{equation}%
We also introduce the fractional functional space $\mathcal{W}^{\lambda
}(Q_{T})$ having the inner product%
\begin{equation}
(\mathcal{Q},\mathcal{Q}^{\ast })_{\mathcal{W}^{\lambda
}(Q_{T})}=\int\limits_{0}^{T}(\mathcal{Q}(.t),\mathcal{Q}^{\ast
}(.,t))_{H_{\ \rho }^{1}(\Omega )}dt+\int\limits_{0}^{T}(^{C}\partial
_{0t}^{\lambda }\mathcal{Q}(.,t),^{C}\partial _{0t}^{\lambda }\mathcal{Q}%
^{\ast }(.,t))_{H_{\ \rho }^{1}(\Omega )}dt,  \label{e2.4*}
\end{equation}%
and with norm%
\begin{equation}
\Vert \mathcal{Q}\Vert _{\mathcal{W}^{\lambda }(Q_{T})}^{2}=\Vert \mathcal{Q}%
\Vert _{L^{2}(0,T;H_{\ \rho }^{1}(\Omega ))}^{2}+\Vert ^{C}\partial
_{0t}^{\lambda }\mathcal{Q}\Vert _{L^{2}(0,T;H_{\ \rho }^{1}(\Omega ))}^{2}.
\label{e2.4**}
\end{equation}%
We denote by $C(0,T;L^{2}(\Omega ))$ the set of all continuous functions $%
V^{\ast }(.,t):[0,T]\rightarrow L^{2}(\Omega )$ with the norm%
\begin{equation}
\left\Vert V^{\ast }\right\Vert _{C(0,T;L^{2}(\Omega ))}^{2}=\underset{0\leq
t\leq T}{\sup }\left\Vert V^{\ast }(.,t)\right\Vert _{L^{2}(\Omega
)}^{2}<\infty .  \label{e2.4***}
\end{equation}%
We recall some definitions of fractional derivatives and fractional integral
[52, 53]. Let $\Gamma (\cdot )$ denote the Gamma function. For any positive
integer $n$ where: $n-1<\alpha <n,$ the Caputo derivative, and fractional
integral of order $\alpha $ are respectively defined by\newline
The left Caputo derivative
\begin{equation}
^{C}\partial _{0t}^{\alpha }v(t)=\frac{1}{\Gamma (n-\alpha )}%
\int\limits_{0}^{t}\frac{v^{(n)}(\tau )}{(t-\tau )^{\alpha -n+1}}d\tau
,\quad \forall t\in \lbrack 0,T],  \label{e2.5}
\end{equation}%
The right Caputo derivative
\begin{equation}
^{C}\partial _{T}^{\alpha }v(t)=\frac{(-1)^{n}}{\Gamma (n-\alpha )}%
\int\limits_{t}^{T}\frac{v^{(n)}(\tau )}{(\tau -t)^{\alpha -n+1}}d\tau
,\quad \forall t\in \lbrack 0,T],  \label{e2.6}
\end{equation}%
and the fractional integral
\begin{equation}
I_{t}^{\alpha }v(t)=D_{0t}^{-\alpha }v(t)=\frac{1}{\Gamma (\alpha )}%
\int\limits_{t}^{T}\frac{v(\tau )}{(t-\tau )^{1-\alpha }}d\tau ,\quad
\forall t\in \lbrack 0,T].  \label{e2.9}
\end{equation}

\begin{lem}
\lbrack 54]. Let a nonnegative absolutely continuous function $\mathcal{P}%
(t) $ satisfy the inequality
\begin{equation*}
^{C}\partial _{0t}^{\beta }\mathcal{P}(t)\leq C\mathcal{P}%
(t)+k(t),~~~0<\beta <1,
\end{equation*}%
for almost all $t\in \lbrack 0,T]$, where $C$ is positive and $k(t)$ is an
integrable nonnegative function on $[0,T]$. Then
\begin{equation*}
\mathcal{P}(t)\leq \mathcal{P}(0)E_{\beta }(Ct^{\beta })+\Gamma (\beta
)E_{\beta ,\beta }(Ct^{\beta })D_{0t}^{-\beta }k(t),
\end{equation*}%
where
\begin{equation*}
E_{\beta }(x)=\sum_{n=0}^{\infty }\frac{x^{n}}{\Gamma (\beta n+1)},\text{
and }E_{\beta ,\alpha }(x)=\sum_{n=0}^{\infty }\frac{x^{n}}{\Gamma (\beta
n+\alpha )},
\end{equation*}%
are the Mittag-Leffler functions.
\end{lem}

\begin{lem}
\lbrack 54]. For any absolutely continuous function $v(t)$ on $[0,T]$, the
following inequality holds
\begin{equation*}
v(t)^{C}\partial _{0t}^{\alpha }v(t)\geq \frac{1}{2}~^{C}\partial
_{0t}^{\alpha }v^{2}(t),~~~0<\alpha <1.
\end{equation*}
\end{lem}

We use the following Gronwall-Bellman lemma

\begin{lem}
\lbrack 55] Let $R(s)$ be nonnegative and absolutely continuous on $[0,T]$,
and suppose that for almost all $s\in \lbrack 0,T]$, the function $R$
satisfies the inequality
\begin{equation}
\frac{dR}{ds}\leq J(s)R(s)+I(s),  \label{e2.12}
\end{equation}%
where the functions $J(s)$ and $I(s)$ are summable and nonnegative on $%
[0,T]. $ Then%
\begin{equation}
R(s)\leq \exp \left\{ \int\limits_{0}^{s}J(t)dt\right\} \left(
R(0)+\int\limits_{0}^{s}I(t)dt\right) .  \label{e2.13}
\end{equation}
\end{lem}

We also use the following inequality [54]
\begin{equation}
D_{0t}^{-\alpha }\Vert f\Vert _{L_{\ \rho }^{2}(\Omega )}^{2}\leq \dfrac{%
t^{\alpha -1}}{\Gamma (\alpha )}\int\limits_{0}^{t}\Vert f\Vert _{L_{\ \rho
}^{2}(\Omega )}^{2}d\tau ,  \label{e2.14}
\end{equation}%
the Cauchy$\varepsilon $-inequality
\begin{equation}
ab\leq \frac{\varepsilon }{2}a^{2}+\frac{1}{2\varepsilon }b^{2},~~~\forall
\varepsilon >0,  \label{e2.15}
\end{equation}%
where $a$ and $b$ are positive numbers. \newline
and the\textbf{\ }Poincare type inequalities [56]
\begin{equation}
\left\Vert \mathcal{J}_{x}(\xi u)\right\Vert _{L^{2}(\Omega )}^{2}\leqslant
\frac{b^{3}}{2}\Vert u(.,t)\Vert _{L_{\rho }^{2}(\Omega )}^{2},
\label{e2.16}
\end{equation}%
\begin{equation}
\left\Vert \mathcal{J}_{x}^{2}(\xi u)\right\Vert _{L^{2}(\Omega
)}^{2}\leqslant \frac{b^{2}}{2}\left\Vert \mathcal{J}_{x}(\xi u)\right\Vert
_{L^{2}(\Omega )}^{2}\leq \frac{b^{5}}{4}\Vert u(.,t)\Vert _{L_{\rho
}^{2}(\Omega )}^{2},  \label{e2.17}
\end{equation}%
where
\begin{equation*}
\mathcal{J}_{x}(\xi v)=\int\limits_{0}^{x}\xi v(\xi ,t)d\xi ,~~\mathcal{J}%
_{x}^{2}(\xi v)=\int\limits_{0}^{x}\int\limits_{0}^{\xi }\eta v(\eta
,t)d\eta .
\end{equation*}

\section{Reformulation of the linear problem}

We consider a fractional coupled system of the form
\begin{equation}
\left\{
\begin{array}{c}
\mathcal{L}_{1}(u,v)=^{C}\partial _{0t}^{\beta }u+-\frac{1}{x}\left(
xu_{x}\right) _{x}-\frac{\partial }{\partial t}\frac{1}{x}\left(
xu_{x}\right) _{x}+z_{1}v+u_{t}=f\left( x,t\right) \\
\mathcal{L}_{2}(u,v)=^{C}\partial _{0t}^{\gamma }v+-\frac{1}{x}\left(
xv_{x}\right) _{x}-\frac{\partial }{\partial t}\frac{1}{x}\left(
xv_{x}\right) _{x}+z_{2}u+v_{t}=g\left( x,t\right)%
\end{array}%
\right.  \label{e3.1}
\end{equation}%
supplemented by the initial conditions
\begin{equation}
\left\{
\begin{array}{c}
\ell _{1}u=u(x,0)=\varphi _{1}(x),~~~\ell _{2}u=u_{t}(x,0)=\varphi _{2}(x)
\\
\ell _{3}v=v(x,0)=\psi _{1}(x),~~~\ell _{4}v=v_{t}(x,0)=\psi _{2}(x),%
\end{array}%
\right.  \label{e3.2}
\end{equation}%
and the Neumann and integral boundary conditions
\begin{equation}
u_{x}(b,t)=0,~v_{x}(b,t)=0,\int\limits_{0}^{b}xudx=0,\text{ }%
\int\limits_{0}^{b}xvdx=0.  \label{e3.3}
\end{equation}%
We assume that there exists a solution $(u,v)\in (C^{2,2}(\overline{Q}))^{2}$
consisting of the set of functions together with their partial derivatives
of order $2$ in $x$ and $t$, which are continuous on $\overline{Q}$.\newline
The solution of system (\ref{e3.1})-(\ref{e3.3}) can be regarded as the
solution of the operator equation $\mathcal{X}W=\mathcal{F},$ where $W,$ $%
\mathcal{X}W$ and $\mathcal{F}$ are respectively the pairs $W=(u,v),$ $%
\mathcal{X}W=\left( L_{1}u,L_{2}v\right) ,$ $\mathcal{F}=\left( \mathcal{F}%
_{1},\mathcal{F}_{2}\right) ,$ with
\begin{equation*}
L_{1}u=\left\{ \mathcal{L}_{1}u,\ell _{1}u,\ell _{2}u\right\} ,\quad
L_{2}v=\left\{ \mathcal{L}_{2}v,\ell _{3}v,\ell _{4}v\right\} ,
\end{equation*}%
and
\begin{equation*}
\mathcal{F}_{1}=\left\{ f,\varphi _{1},\varphi _{2}\right\} ,\quad \mathcal{F%
}_{2}=\left\{ g,\psi _{1},\psi _{2}\right\} .
\end{equation*}%
The operator $\mathcal{X}$ is considered from a space $B$ into a space $H$,
where $B$ is a Banach space consisting of all functions $(u,v)\in \left(
L^{2}(0,T;L_{\rho }^{2}(\Omega ))\right) ^{2}$ satisfying conditions (\ref%
{e3.3}) and having the finite norm
\begin{equation}
\Vert W\Vert _{B}^{2}=\lVert u\rVert _{\mathcal{W}^{\beta
}(Q_{T})}^{2}+\lVert v\rVert _{\mathcal{W}^{\gamma }(Q_{T})}^{2}+\lVert
u\rVert _{C(0,t.,H_{\ \rho }^{1}(\Omega ))}^{2}+\lVert v\rVert _{C(0,t.,H_{\
\rho }^{1}(\Omega ))}^{2},  \notag
\end{equation}%
and $H=\left( L^{2}(Q_{T})\right) ^{2}\times \left( H_{\rho }^{1}(\Omega
)\right) ^{4}$ is the Hilbert space consisting of vector-valued functions $%
\mathcal{S}=\left( \{f,\varphi _{1},\psi _{1}\},\{g,\varphi _{2},\psi
_{2}\}\right) $ with norm
\begin{equation*}
\Vert \mathcal{S}\Vert _{H}^{2}=\Vert f\Vert _{L^{2}(0,T;L_{\rho
}^{2}(\Omega ))}^{2}+\Vert g\Vert _{L^{2}(0,T;L_{\rho }^{2}(\Omega
))}^{2}+\left\Vert \varphi _{1}\right\Vert _{H_{\rho }^{1}(\Omega
)}^{2}+\left\Vert \varphi _{2}\right\Vert _{H_{\rho }^{1}(\Omega
)}^{2}+\left\Vert \psi _{1}\right\Vert _{H_{\rho }^{1}(\Omega
)}^{2}+\left\Vert \psi _{2}\right\Vert _{H_{\rho }^{1}(\Omega )}^{2}.
\end{equation*}%
Let $D(\mathcal{X}),$ be the domain of definition of the operator $\mathcal{X%
},$ defined by:
\begin{equation*}
D(\mathcal{X})=\left\{
\begin{array}{c}
(u,v)\in \left( L^{2}(0,T;L_{\rho }^{2}(\Omega ))\right) ^{2}\text{ such
that }^{C}\partial _{0t}^{\beta }u,^{C}\partial _{0t}^{\gamma }v,u_{x}, \\
v_{x},u_{xx},v_{xx},u_{tx},v_{tx},u_{txx},v_{txx}\in L^{2}(0,T;L_{\rho
}^{2}(\Omega )) \\
u_{x}(b,t)=0,~v_{x}(b,t)=0,\text{ }\int\limits_{0}^{b}xudx=0,\text{ }%
\int\limits_{0}^{b}xvdx=0.%
\end{array}%
\right.
\end{equation*}

\section{Uniqueness of the solution}

In this section, we prove the uniqueness result for the fractional system (%
\ref{e3.1})-(\ref{e3.3}), that is we establish an energy inequality for the
operator $\mathcal{X}$ and we give some of its consequences.

\begin{thm}
For any $(u,v)\in D(\mathcal{X})$, $f,$ $g\in L^{2}(0,T;L_{\ \rho
}^{2}(\Omega )),$ and $\varphi _{1},\psi _{1},\varphi _{2},\psi _{2}\in H_{\
\rho }^{1}(\Omega ),$ the solution of the problem (\ref{e3.1})-(\ref{e3.3})
verifies the a priori bound
\begin{eqnarray}
&&\lVert u\rVert _{\mathcal{W}^{\beta }(Q_{T})}^{2}+\lVert u\rVert _{%
\mathcal{W}^{\gamma }(Q_{T})}^{2}+\lVert u\rVert _{C(0,t.,H_{\ \rho
}^{1}(\Omega )}^{2}+\lVert v\rVert _{C(0,t.,H_{\ \rho }^{1}(\Omega )}^{2}
\notag \\
&\leq &\mathcal{M}\left( \lVert f\rVert _{L^{2}(0,T;L_{\ \rho }^{2}(\Omega
))}^{2}+\lVert g\rVert _{L^{2}(0,T;L_{\ \rho }^{2}(\Omega ))}^{2}+\lVert
\varphi _{1}\rVert _{H_{\ \rho }^{1}(\Omega )}^{2}+\lVert \psi _{1}\rVert
_{H_{\ \rho }^{1}(\Omega )}^{2}\right.  \notag \\
&&+\left. \lVert \varphi _{2}\rVert _{H_{\ \rho }^{1}(\Omega )}^{2}+\lVert
\psi _{2}\rVert _{H_{\ \rho }^{1}(\Omega )}^{2}\right) ,  \label{e4.1}
\end{eqnarray}%
where $\mathcal{M=Y}^{\ast \ast }e^{T\mathcal{Y}^{\ast \ast }}$ is a positve
constant with%
\begin{equation}
\left\{
\begin{array}{c}
\mathcal{Y}^{\ast \ast }=\max \left( 1,\mathcal{Y}^{\ast }\right) ,\text{ }%
\mathcal{Y}^{\ast }=\frac{\mathcal{Y}}{\min \left( 1,\frac{b^{2}}{4}\right) }%
,\text{ }\mathcal{Y=}\chi ^{\ast }\chi \max \left( \frac{T^{\beta -1}}{%
\Gamma (\beta )},\frac{T^{\gamma -1}}{\Gamma (\gamma )}\right) \\
\chi ^{\ast }=\Gamma (\beta -1)E_{\beta -1,\beta -1}(\chi t^{\beta -1})\max
\left\{ 1,\frac{T^{\beta -1}}{(\beta -1)\Gamma (\beta -1)}\right\} ,\text{ }%
\chi =D^{\ast \ast }\left( 1+D^{\ast \ast }e^{D^{\ast \ast }T}\right) \\
D^{\ast \ast }=D^{\ast }\max \left\{ 1,\frac{b^{4}}{2},\frac{T^{2-\beta }}{%
(2-\beta )\Gamma (2-\beta )},\frac{T^{2-\gamma }}{(2-\gamma )\Gamma
(2-\gamma )}\right\} \\
D^{\ast }=2\max \left\{ 3,\frac{b^{6}}{8}+\frac{1}{2},\frac{b^{4}}{8}+\frac{5%
}{2}\right\} .%
\end{array}%
\right.  \label{e4.2}
\end{equation}
\end{thm}

\textbf{Proof. }The fractional partial differential equations in (\ref{e3.1}%
), and the following fractional integro-differential operators
\begin{equation*}
\mathcal{M}_{1}u=^{C}\partial _{0t}^{\beta }u+u_{t}-\mathcal{J}_{x}^{2}(\xi
u_{t}),\text{ and }\mathcal{M}_{2}v=^{C}\partial _{0t}^{\gamma }v+v_{t}-%
\mathcal{J}_{x}^{2}(\xi v_{t}),
\end{equation*}%
lead to
\begin{eqnarray}
&&2\left( ^{C}\partial _{0t}^{\beta }u,u_{t}\right) _{{L_{\rho \ }^{2}\left(
\Omega \right) }}-\left( ^{C}\partial _{0t}^{\beta }u,\mathcal{J}%
_{x}^{2}(\xi u_{t})\right) _{{L_{\ \rho }^{2}\left( \Omega \right) }}-\left(
\dfrac{1}{x}\left( xu_{x}\right) _{x},u_{t}\right) _{{L_{\ \rho }^{2}\left(
\Omega \right) }}+\left( \dfrac{1}{x}\left( xu_{x}\right) _{x},\mathcal{J}%
_{x}^{2}(\xi u_{t})\right) _{{L_{\ \rho }^{2}\left( \Omega \right) }}  \notag
\\
&&+\left( ^{C}\partial _{0t}^{\beta }u,^{C}\partial _{0t}^{\beta }u\right) _{%
{L_{\ \rho }^{2}\left( \Omega \right) }}-\left( \dfrac{1}{x}\left(
xu_{x}\right) _{x},^{C}\partial _{0t}^{\beta }u\right) _{{L_{\ \rho
}^{2}\left( \Omega \right) }}-\left( \frac{1}{x}\left( xu_{x}\right)
_{xt},^{C}\partial _{0t}^{\beta }u\right) _{{L_{\ \rho }^{2}\left( \Omega
\right) }}  \notag \\
&&+\left( ^{C}\partial _{0t}^{\beta }u,z_{1}v\right) _{{L_{\ \rho
}^{2}\left( \Omega \right) }}+\left( ^{C}\partial _{0t}^{\gamma
}v,^{C}\partial _{0t}^{\gamma }v\right) _{{L_{\ \rho }^{2}\left( \Omega
\right) }}-\left( \dfrac{1}{x}\left( xv_{x}\right) _{x},^{C}\partial
_{0t}^{\gamma }v\right) _{{L_{\ \rho }^{2}\left( \Omega \right) }}-\left(
\frac{1}{x}\left( xv_{x}\right) _{xt},^{C}\partial _{0t}^{\gamma }v\right) _{%
{L_{\ \rho }^{2}\left( \Omega \right) }}  \notag \\
&&+2\left( ^{C}\partial _{0t}^{\gamma }v,v_{t}\right) _{{L_{\ \rho
}^{2}\left( \Omega \right) }}+\left( ^{C}\partial _{0t}^{\gamma
}v,z_{2}u\right) _{{L_{\ \rho }^{2}\left( \Omega \right) }}-\left( \dfrac{1}{%
x}\left( xu_{x}\right) _{xt},u_{t}\right) _{{L_{\ \rho }^{2}\left( \Omega
\right) }}+\left( \dfrac{1}{x}\left( xu_{x}\right) _{xt},\mathcal{J}%
_{x}^{2}(\xi u_{t})\right) _{{L_{\ \rho }^{2}\left( \Omega \right) }}  \notag
\\
&&+\left( z_{1}v,u_{t}\right) _{L_{\ \rho }^{2}\left( \Omega \right)
}-\left( z_{1}v,\mathcal{J}_{x}^{2}(\xi u_{t})\right) _{L_{\ \rho
}^{2}\left( \Omega \right) }-\left( ^{C}\partial _{0t}^{\gamma }v,\mathcal{J}%
_{x}^{2}(\xi v_{t})\right) _{{L_{\ \rho }^{2}\left( \Omega \right) }}  \notag
\\
&&-\left( \dfrac{1}{x}\left( xv_{x}\right) _{x},v_{t}\right) _{{L_{\ \rho
}^{2}\left( \Omega \right) }}+\left( \dfrac{1}{x}\left( xv_{x}\right) _{x},%
\mathcal{J}_{x}^{2}(\xi v_{t})\right) _{{L_{\ \rho }^{2}\left( \Omega
\right) }}-\left( \dfrac{1}{x}\left( xv_{x}\right) _{xt},v_{t}\right) _{{%
L_{\ \rho }^{2}\left( \Omega \right) }}  \notag \\
&&+\left( \dfrac{1}{x}\left( xv_{x}\right) _{xt},\mathcal{J}_{x}^{2}(\xi
v_{t})\right) _{{L_{\ \rho }^{2}\left( \Omega \right) }}+\left(
z_{2}u,v_{t}\right) _{L_{\ \rho }^{2}\left( \Omega \right) }-\left( z_{2}u,%
\mathcal{J}_{x}^{2}(\xi v_{t})\right) _{L_{\ \rho }^{2}\left( \Omega \right)
}  \notag \\
&&+\left\Vert u_{t}\right\Vert _{L_{\ \rho }^{2}\left( \Omega \right)
}^{2}+\left\Vert v_{t}\right\Vert _{L_{\ \rho }^{2}\left( \Omega \right)
}^{2}-\left( u_{t},\mathcal{J}_{x}^{2}(\xi u_{t})\right) _{L_{\ \rho
}^{2}\left( \Omega \right) }-\left( v_{t},\mathcal{J}_{x}^{2}(\xi
v_{t})\right) _{L_{\ \rho }^{2}\left( \Omega \right) }  \notag \\
&=&\left( f,u_{t}\right) _{L_{\rho }^{2}\left( \Omega \right) }-\left( f,%
\mathcal{J}_{x}^{2}(\xi u_{t})\right) _{L_{\rho }^{2}\left( \Omega \right)
}+(g,v_{t})_{L_{\rho }^{2}\left( \Omega \right) }-(g,\mathcal{J}_{x}^{2}(\xi
v_{t}))_{L_{\rho }^{2}\left( \Omega \right) }+\left( f,^{C}\partial
_{0t}^{\beta }u\right) _{L_{\rho }^{2}\left( \Omega \right) }  \label{e4.2*}
\\
&&+\left( g,^{C}\partial _{0t}^{\gamma }v\right) _{L_{\rho }^{2}\left(
\Omega \right) }.  \notag
\end{eqnarray}%
Using boundary conditions (\ref{e3.3}), we evaluate the following terms on
the LHS of (\ref{e4.2*}) as follows
\begin{eqnarray}
-\left( ^{C}\partial _{0t}^{\beta }u,\mathcal{J}_{x}^{2}(\xi u_{t})\right) _{%
{L_{\ \rho }^{2}\left( \Omega \right) }} &=&\left( ^{C}\partial _{0t}^{\beta
}\left( \mathcal{J}_{x}(\xi u)\right) ,\mathcal{J}_{x}(\xi u_{t})\right) _{{%
L^{2}\left( \Omega \right) }},  \label{e4.3} \\
-\left( \dfrac{1}{x}\left( xu_{x}\right) _{x},u_{t}\right) _{{L_{\ \rho
}^{2}\left( \Omega \right) }} &=&\dfrac{1}{2}\dfrac{\partial }{\partial t}%
\lVert u_{x}\rVert _{L_{\ \rho }^{2}(\Omega )}^{2},  \label{e4.4} \\
\left( \dfrac{1}{x}\left( xu_{x}\right) _{x},\mathcal{J}_{x}^{2}(\xi
u_{t})\right) _{{L_{\ \rho }^{2}\left( \Omega \right) }} &=&-\left( u_{x},%
\mathcal{J}_{x}(\xi u_{t})\right) _{{L_{\ \rho }^{2}\left( \Omega \right) }},
\label{e4.5} \\
-\left( \dfrac{1}{x}\left( xu_{x}\right) _{xt},u_{t}\right) _{{L_{\ \rho
}^{2}\left( \Omega \right) }} &=&\lVert u_{xt}\rVert _{L_{\ \rho
}^{2}(\Omega )}^{2},  \label{e4.6} \\
\left( \dfrac{1}{x}\left( xu_{x}\right) _{xt},\mathcal{J}_{x}^{2}(\xi
u_{t})\right) _{{L_{\ \rho }^{2}\left( \Omega \right) }} &=&-\left( u_{xt},%
\mathcal{J}_{x}(\xi u_{t})\right) _{{L_{\ \rho }^{2}\left( \Omega \right) }},
\label{e4.7} \\
-\left( z_{1}v,\mathcal{J}_{x}^{2}(\xi u_{t})\right) _{L_{\ \rho }^{2}\left(
\Omega \right) } &=&-z_{1}\left( \mathcal{J}_{x}^{2}(\xi v),u_{t}\right)
_{L_{\ \rho }^{2}\left( \Omega \right) },  \label{e4.8}
\end{eqnarray}%
\begin{equation}
\left( ^{C}\partial _{0t}^{\beta }u,^{C}\partial _{0t}^{\beta }u\right) _{{%
L_{\ \rho }^{2}\left( \Omega \right) }}=\lVert ^{C}\partial _{0t}^{\beta
-1}u_{t}\rVert _{L_{\ \rho }^{2}(\Omega )}^{2},  \label{e4.8*}
\end{equation}%
\begin{equation}
-\left( u_{t},\mathcal{J}_{x}^{2}(\xi u_{t})\right) _{L_{\ \rho }^{2}\left(
\Omega \right) }=\left\Vert \mathcal{J}_{x}(\xi u_{t})\right\Vert _{L_{\
}^{2}\left( \Omega \right) }^{2},  \label{e4.8**}
\end{equation}%
\begin{equation}
-\left( v_{t},\mathcal{J}_{x}^{2}(\xi v_{t})\right) _{L_{\ \rho }^{2}\left(
\Omega \right) }=\left\Vert \mathcal{J}_{x}(\xi v_{t})\right\Vert _{L_{\
}^{2}\left( \Omega \right) }^{2},  \label{e4.8***}
\end{equation}%
\begin{equation}
-\left( \dfrac{1}{x}\left( xu_{x}\right) _{x},^{C}\partial _{0t}^{\beta
}u\right) _{{L_{\ \rho }^{2}\left( \Omega \right) }}=\left( ^{C}\partial
_{0t}^{\beta -1}u_{t},u_{x}\right) _{{L_{\ \rho }^{2}\left( \Omega \right) }%
},  \label{e4.9*}
\end{equation}%
\begin{equation}
-\left( \frac{1}{x}\left( xu_{x}\right) _{xt},^{C}\partial _{0t}^{\beta
}u\right) _{{L_{\ \rho }^{2}\left( \Omega \right) }}=\left( ^{C}\partial
_{0t}^{\beta -1}u_{xt},u_{xt}\right) _{{L_{\ \rho }^{2}\left( \Omega \right)
}},  \label{e4.11*}
\end{equation}%
\begin{equation}
\left( ^{C}\partial _{0t}^{\beta }u,z_{1}v\right) _{{L_{\ \rho }^{2}\left(
\Omega \right) }}=\left( ^{C}\partial _{0t}^{\beta -1}u_{t},z_{1}v\right) _{{%
L_{\ \rho }^{2}\left( \Omega \right) }}.  \label{e4.13*}
\end{equation}%
In the same fashion, we have the equations (\ref{e4.3})-(\ref{e4.13*}) with $%
\beta $ replaced by $\gamma $, and $u$ replaced by $v.$ Since $0$ $<\beta
-1<1,$ then by using Lemma 2.2, we have
\begin{equation}
2\left( ^{C}\partial _{0t}^{\beta }u,u_{t}\right) _{{L_{\ \rho }^{2}\left(
\Omega \right) }}=2\left( ^{C}\partial _{0t}^{\beta -1}u_{t},u_{t}\right) _{{%
L_{\ \rho }^{2}\left( \Omega \right) }}\geq ^{C}\partial _{0t}^{\beta
-1}\lVert u_{t}\rVert _{L_{\ \rho }^{2}(\Omega )}^{2},  \label{e4.9}
\end{equation}%
\begin{eqnarray}
\left( ^{C}\partial _{0t}^{\beta }\left( \mathcal{J}_{x}(\xi u)\right) ,%
\mathcal{J}_{x}(\xi u_{t})\right) _{{L^{2}\left( \Omega \right) }} &=&\left(
^{C}\partial _{0t}^{\beta -1}\left( \mathcal{J}_{x}(\xi u_{t})\right) ,%
\mathcal{J}_{x}(\xi u_{t})\right) _{{L^{2}\left( \Omega \right) }}  \notag \\
&\geq &\dfrac{1}{2}~~^{C}\partial _{0t}^{\beta -1}\lVert \mathcal{J}_{x}(\xi
u_{t})\rVert _{L_{\ \rho }^{2}(\Omega )}^{2}.  \label{e4.10}
\end{eqnarray}%
\begin{equation}
\left( ^{C}\partial _{0t}^{\beta -1}u_{xt},u_{xt}\right) _{{L_{\ \rho
}^{2}\left( \Omega \right) }}\geq \frac{1}{2}^{C}\partial _{0t}^{\beta
-1}\lVert u_{xt}\rVert _{L_{\ \rho }^{2}(\Omega )}^{2},  \label{e4.10*}
\end{equation}%
Combination of (\ref{e4.2})-(\ref{e4.10*}) yields%
\begin{eqnarray}
&&\left\Vert ^{C}\partial _{0t}^{\beta -1}u_{t}\right\Vert _{L_{\ \rho
}^{2}(\Omega )}^{2}+\left\Vert ^{C}\partial _{0t}^{\gamma
-1}v_{t}\right\Vert _{L_{\ \rho }^{2}(\Omega )}^{2}+\frac{1}{2}%
~~^{C}\partial _{0t}^{\beta -1}\lVert u_{t}\rVert _{L_{\ \rho }^{2}(\Omega
)}^{2}+\frac{1}{2}~^{C}\partial _{0t}^{\gamma -1}\lVert v_{t}\rVert _{L_{\
\rho }^{2}(\Omega )}^{2}  \notag \\
&&+~\frac{1}{2}~~^{C}\partial _{0t}^{\beta -1}\lVert \mathcal{J}_{x}(\xi
u_{t})\rVert _{L_{\ \rho }^{2}(\Omega )}^{2}+\frac{1}{2}~^{C}\partial
_{0t}^{\gamma -1}\lVert \mathcal{J}_{x}(\xi v_{t})\rVert _{L_{\ \rho
}^{2}(\Omega )}^{2}+\frac{1}{2}\dfrac{\partial }{\partial t}\lVert
u_{x}\rVert _{L_{\ \rho }^{2}(\Omega )}^{2}+\frac{1}{2}\dfrac{\partial }{%
\partial t}\lVert v_{x}\rVert _{L_{\ \rho }^{2}(\Omega )}^{2}  \notag \\
&&+\dfrac{1}{2}~~^{C}\partial _{0t}^{\beta -1}\lVert u_{xt}\rVert _{L_{\
\rho }^{2}(\Omega )}^{2}+\dfrac{1}{2}~~^{C}\partial _{0t}^{\gamma -1}\lVert
v_{xt}\rVert _{L_{\ \rho }^{2}(\Omega )}^{2}+\left\Vert u_{t}\right\Vert
_{L_{\ \rho }^{2}\left( \Omega \right) }^{2}+\left\Vert v_{t}\right\Vert
_{L_{\ \rho }^{2}\left( \Omega \right) }^{2}  \notag \\
&&+\left\Vert \mathcal{J}_{x}(\xi u_{t})\right\Vert _{L_{\ }^{2}\left(
\Omega \right) }^{2}+\left\Vert \mathcal{J}_{x}(\xi v_{t})\right\Vert _{L_{\
}^{2}\left( \Omega \right) }^{2}  \notag \\
&\leq &\left( f,u_{t}\right) _{L_{\rho }^{2}\left( \Omega \right) }-\left( f,%
\mathcal{J}_{x}^{2}(\xi u_{t})\right) _{L_{\rho }^{2}\left( \Omega \right)
}+(g,v_{t})_{L_{\rho }^{2}\left( \Omega \right) }-(g,\mathcal{J}_{x}^{2}(\xi
v_{t}))_{L_{\rho }^{2}\left( \Omega \right) }+\left( f,^{C}\partial
_{0t}^{\beta -1}u_{t}\right) _{L_{\rho }^{2}\left( \Omega \right) }  \notag
\\
&&+\left( g,^{C}\partial _{0t}^{\gamma -1}v_{t}\right) _{L_{\rho }^{2}\left(
\Omega \right) }-\left( ^{C}\partial _{0t}^{\beta -1}u_{t},z_{1}v\right) _{{%
L_{\ \rho }^{2}\left( \Omega \right) }}-\left( ^{C}\partial _{0t}^{\gamma
-1}v_{t},z_{2}u\right) _{{L_{\ \rho }^{2}\left( \Omega \right) }}  \notag \\
&&-\left( ^{C}\partial _{0t}^{\beta -1}u_{t},u_{x}\right) _{{L_{\ \rho
}^{2}\left( \Omega \right) }}-\left( ^{C}\partial _{0t}^{\gamma
-1}v_{t},v_{x}\right) _{{L_{\ \rho }^{2}\left( \Omega \right) }}  \notag \\
&&-\left( z_{1}v,u_{t}\right) _{L_{\ \rho }^{2}\left( \Omega \right)
}+\left( u_{x},\mathcal{J}_{x}(\xi u_{t})\right) _{{L_{\ \rho }^{2}\left(
\Omega \right) }}+\left( u_{xt},\mathcal{J}_{x}(\xi u_{t})\right) _{{L_{\
\rho }^{2}\left( \Omega \right) }}+z_{1}\left( \mathcal{J}_{x}^{2}(\xi
v),u_{t}\right) _{L_{\ \rho }^{2}\left( \Omega \right) }  \notag \\
&&-\left( z_{2}u,v_{t}\right) _{L_{\ \rho }^{2}\left( \Omega \right)
}+\left( v_{x},\mathcal{J}_{x}(\xi v_{t})\right) _{{L_{\ \rho }^{2}\left(
\Omega \right) }}+\left( v_{xt},\mathcal{J}_{x}(\xi v_{t})\right) _{{L_{\
\rho }^{2}\left( \Omega \right) }}+z_{2}\left( \mathcal{J}_{x}^{2}(\xi
u),v_{t}\right) _{L_{\ \rho }^{2}\left( \Omega \right) }.  \label{e4.11}
\end{eqnarray}%
By applying Cauchy-$\varepsilon $-inequality (\ref{e2.15}) and Poincare type
inequalities (\ref{e2.16}) and (\ref{e2.17}) to the the right-hand side of (%
\ref{e4.11}), we obtain the inequalities%
\begin{eqnarray}
\left( f,u_{t}\right) _{L_{\rho }^{2}\left( \Omega \right) } &\leq &\dfrac{%
\eta _{1}}{2}\lVert f\rVert _{L_{\ \rho }^{2}(\Omega )}^{2}+\dfrac{1}{2\eta
_{1}}\lVert u_{t}\rVert _{L_{\ \rho }^{2}(\Omega )}^{2},  \label{e4.12} \\
-\left( f,\mathcal{J}_{x}^{2}(\xi u_{t})\right) _{L_{\rho }^{2}\left( \Omega
\right) } &\leq &\dfrac{1}{2\eta _{2}}\lVert f\rVert _{L_{\ \rho
}^{2}(\Omega )}^{2}+\dfrac{\eta _{2}b^{6}}{8}\lVert u_{t}\rVert _{L_{\ \rho
}^{2}(\Omega )}^{2},  \label{e4.13} \\
-z_{1}\left( v,u_{t}\right) _{L_{\ \rho }^{2}\left( \Omega \right) } &\leq &%
\dfrac{z_{1}^{2}}{2\eta _{3}}\lVert v\rVert _{L_{\ \rho }^{2}(\Omega )}^{2}+%
\dfrac{\eta _{3}}{2}\lVert u_{t}\rVert _{L_{\ \rho }^{2}(\Omega )}^{2},
\label{e4.14} \\
\left( u_{x},\mathcal{J}_{x}(\xi u_{t})\right) _{{L_{\ \rho }^{2}\left(
\Omega \right) }} &\leq &\dfrac{1}{2}\lVert u_{x}\rVert _{L_{\ \rho
}^{2}(\Omega )}^{2}+\dfrac{1}{2}\lVert \mathcal{J}_{x}(\xi u_{t})\rVert
_{L_{\ \rho }^{2}(\Omega )}^{2},  \label{e4.15} \\
\left( u_{xt},\mathcal{J}_{x}(\xi u_{t})\right) _{{L_{\ \rho }^{2}\left(
\Omega \right) }} &\leq &\dfrac{\eta _{4}}{2}\lVert u_{xt}\rVert _{L_{\ \rho
}^{2}(\Omega )}^{2}+\dfrac{1}{2\eta _{4}}\lVert \mathcal{J}_{x}(\xi
u_{t})\rVert _{L_{\ \rho }^{2}(\Omega )}^{2},  \label{e4.16} \\
z_{1}\left( \mathcal{J}_{x}^{2}(\xi v),u_{t}\right) _{L_{\ \rho }^{2}\left(
\Omega \right) } &\leq &\dfrac{z_{1}^{2}b^{4}}{8\eta _{5}}\lVert v\rVert
_{L_{\ \rho }^{2}(\Omega )}^{2}+\dfrac{\eta _{5}}{2}\lVert u_{t}\rVert
_{L_{\ \rho }^{2}(\Omega )}^{2},  \label{e4.17}
\end{eqnarray}%
\begin{equation}
\left( f,^{C}\partial _{0t}^{\beta -1}u_{t}\right) _{L_{\rho }^{2}\left(
\Omega \right) }\leq \dfrac{\eta _{11}}{2}\lVert ^{C}\partial _{0t}^{\beta
-1}u_{t}\rVert _{L_{\ \rho }^{2}(\Omega )}^{2}+\dfrac{1}{2\eta _{11}}\lVert
f\rVert _{L_{\ \rho }^{2}(\Omega )}^{2},  \label{e4.17**}
\end{equation}%
\begin{equation}
\left( g,^{C}\partial _{0t}^{\gamma -1}v_{t}\right) _{L_{\rho }^{2}\left(
\Omega \right) }\leq \dfrac{\eta _{12}}{2}\lVert ^{C}\partial _{0t}^{\gamma
-1}v_{t}\rVert _{L_{\ \rho }^{2}(\Omega )}^{2}+\dfrac{1}{2\eta _{12}}\lVert
g\rVert _{L_{\ \rho }^{2}(\Omega )}^{2},  \label{e4.18*}
\end{equation}%
\begin{equation}
-z_{1}\left( ^{C}\partial _{0t}^{\beta -1}u_{t},v\right) _{{L_{\ \rho
}^{2}\left( \Omega \right) }}\leq \dfrac{\eta _{13}}{2}\lVert ^{C}\partial
_{0t}^{\beta -1}u_{t}\rVert _{L_{\ \rho }^{2}(\Omega )}^{2}+\dfrac{z_{1}^{2}%
}{2\eta _{13}}\lVert v\rVert _{L_{\ \rho }^{2}(\Omega )}^{2},  \label{e4.19*}
\end{equation}%
\begin{equation}
-\left( ^{C}\partial _{0t}^{\gamma -1}v_{t},u\right) _{{L_{\ \rho
}^{2}\left( \Omega \right) }}\leq \dfrac{\eta _{14}}{2}\lVert ^{C}\partial
_{0t}^{\gamma -1}v_{t}\rVert _{L_{\ \rho }^{2}(\Omega )}^{2}+\dfrac{1}{2\eta
_{14}}\lVert u\rVert _{L_{\ \rho }^{2}(\Omega )}^{2},  \label{e4.20*}
\end{equation}%
\begin{equation}
-\left( ^{C}\partial _{0t}^{\beta -1}u_{t},u_{x}\right) _{{L_{\ \rho
}^{2}\left( \Omega \right) }}\leq \dfrac{\eta _{15}}{2}\lVert ^{C}\partial
_{0t}^{\beta -1}u_{t}\rVert _{L_{\ \rho }^{2}(\Omega )}^{2}+\dfrac{1}{2\eta
_{15}}\lVert u_{x}\rVert _{L_{\ \rho }^{2}(\Omega )}^{2},  \label{e4.21*}
\end{equation}%
\begin{equation}
-z_{2}\left( ^{C}\partial _{0t}^{\gamma -1}v_{t},v_{x}\right) _{{L_{\ \rho
}^{2}\left( \Omega \right) }}\leq \dfrac{\eta _{16}}{2}\lVert ^{C}\partial
_{0t}^{\gamma -1}v_{t}\rVert _{L_{\ \rho }^{2}(\Omega )}^{2}+\dfrac{z_{2}^{2}%
}{2\eta _{15}}\lVert v_{x}\rVert _{L_{\ \rho }^{2}(\Omega )}^{2},
\label{e4.22*}
\end{equation}%
\begin{equation}
(g,v_{t})_{L_{\rho }^{2}\left( \Omega \right) }\leq \dfrac{\eta _{6}}{2}%
\lVert g\rVert _{L_{\ \rho }^{2}(\Omega )}^{2}+\dfrac{1}{2\eta _{6}}\lVert
v_{t}\rVert _{L_{\ \rho }^{2}(\Omega )}^{2},  \label{e4.12*}
\end{equation}%
\begin{equation}
-\left( g,\mathcal{J}_{x}^{2}(\xi v_{t})\right) _{L_{\rho }^{2}\left( \Omega
\right) }\leq \dfrac{1}{2\eta _{7}}\lVert g\rVert _{L_{\ \rho }^{2}(\Omega
)}^{2}+\dfrac{\eta _{7}b^{6}}{8}\lVert v_{t}\rVert _{L_{\ \rho }^{2}(\Omega
)}^{2},  \label{e4.13**}
\end{equation}%
\begin{equation}
-z_{2}\left( u,v_{t}\right) _{L_{\ \rho }^{2}\left( \Omega \right) }\leq
\dfrac{z_{2}^{2}}{2\eta _{8}}\lVert u\rVert _{L_{\ \rho }^{2}(\Omega )}^{2}+%
\dfrac{\eta _{8}}{2}\lVert v_{t}\rVert _{L_{\ \rho }^{2}(\Omega )}^{2},
\label{e4.14*}
\end{equation}%
\begin{equation}
\left( v_{x},\mathcal{J}_{x}(\xi v_{t})\right) _{{L_{\ \rho }^{2}\left(
\Omega \right) }}\leq \dfrac{1}{2}\lVert v_{x}\rVert _{L_{\ \rho
}^{2}(\Omega )}^{2}+\dfrac{1}{2}\lVert \mathcal{J}_{x}(\xi v_{t})\rVert
_{L_{\ \rho }^{2}(\Omega )}^{2},  \label{e4.15*}
\end{equation}%
\begin{equation}
\left( v_{xt},\mathcal{J}_{x}(\xi v_{t})\right) _{{L_{\ \rho }^{2}\left(
\Omega \right) }}\leq \dfrac{\eta _{9}}{2}\lVert v_{xt}\rVert _{L_{\ \rho
}^{2}(\Omega )}^{2}+\dfrac{1}{2\eta _{9}}\lVert \mathcal{J}_{x}(\xi
v_{t})\rVert _{L_{\ \rho }^{2}(\Omega )}^{2},  \label{e4.16*}
\end{equation}%
\begin{equation}
z_{2}\left( \mathcal{J}_{x}^{2}(\xi u),v_{t}\right) _{L_{\ \rho }^{2}\left(
\Omega \right) }\leq \dfrac{z_{2}^{2}b^{4}}{8\eta _{10}}\lVert u\rVert
_{L_{\ \rho }^{2}(\Omega )}^{2}+\dfrac{\eta _{10}}{2}\lVert v_{t}\rVert
_{L_{\ \rho }^{2}(\Omega )}^{2}.  \label{e4.17*}
\end{equation}%
By inserting (\ref{e4.12})-(\ref{e4.17*}) into (\ref{e4.11}), and taking $%
\eta _{1}=\eta _{2}=\eta _{3}=\eta _{5}=\eta _{6}=\eta _{7}=\eta _{8}=\eta
_{10}=1,$ $\eta _{4}=\eta _{9}=1,$ $\eta _{11}=\eta _{12}=\eta _{13}=\eta
_{14}=\eta _{15}=\eta _{16}=1/4,$ gives
\begin{eqnarray}
&&~\lVert ^{C}\partial _{0t}^{\beta -1}u_{t}\rVert _{L_{\ \rho }^{2}(\Omega
)}^{2}+\lVert ^{C}\partial _{0t}^{\gamma -1}v_{t}\rVert _{L_{\ \rho
}^{2}(\Omega )}^{2}+^{C}\partial _{0t}^{\beta -1}\lVert u_{t}\rVert _{L_{\
\rho }^{2}(\Omega )}^{2}  \notag \\
&&+^{C}\partial _{0t}^{\gamma -1}\lVert v_{t}\rVert _{L_{\ \rho }^{2}(\Omega
)}^{2}+\dfrac{\partial }{\partial t}\lVert u_{x}\rVert _{L_{\ \rho
}^{2}(\Omega )}^{2}+\dfrac{\partial }{\partial t}\lVert v_{x}\rVert _{L_{\
\rho }^{2}(\Omega )}^{2}  \notag \\
&&+^{C}\partial _{0t}^{\beta -1}\lVert \mathcal{J}_{x}(\xi u_{t})\rVert
_{L_{\ \rho }^{2}(\Omega )}^{2}+^{C}\partial _{0t}^{\gamma -1}\lVert
\mathcal{J}_{x}(\xi v_{t})\rVert _{L_{\ \rho }^{2}(\Omega )}^{2}  \notag \\
&&+^{C}\partial _{0t}^{\beta -1}\lVert u_{tx})\rVert _{L_{\ \rho
}^{2}(\Omega )}^{2}+^{C}\partial _{0t}^{\gamma -1}\lVert v_{tx})\rVert
_{L_{\ \rho }^{2}(\Omega )}^{2}  \notag \\
&\leq &D^{\ast }\left( \lVert u_{t}\rVert _{L_{\ \rho }^{2}(\Omega
)}^{2}+\lVert v_{t}\rVert _{L_{\ \rho }^{2}(\Omega )}^{2}+\lVert u\rVert
_{L_{\ \rho }^{2}(\Omega )}^{2}+\lVert v\rVert _{L_{\ \rho }^{2}(\Omega
)}^{2}+\lVert u_{x}\rVert _{L_{\ \rho }^{2}(\Omega )}^{2}\right.  \notag \\
&&+\lVert v_{x}\rVert _{L_{\ \rho }^{2}(\Omega )}^{2}+\lVert u_{xt}\rVert
_{L_{\ \rho }^{2}(\Omega )}^{2}+\lVert v_{xt}\rVert _{L_{\ \rho }^{2}(\Omega
)}^{2}+\lVert \mathcal{J}_{x}(\xi u_{t})\rVert _{L_{\ \rho }^{2}(\Omega
)}^{2}  \notag \\
&&+\left. \lVert \mathcal{J}_{x}(\xi v_{t})\rVert _{L_{\ \rho }^{2}(\Omega
)}^{2}+\lVert f\rVert _{L_{\ \rho }^{2}(\Omega )}^{2}+\lVert g\rVert _{L_{\
\rho }^{2}(\Omega )}^{2}\right) ,  \label{e4.18}
\end{eqnarray}%
where
\begin{equation}
D^{\ast }=2\max \left\{ 3,\frac{b^{6}}{8}+\frac{3}{2},\frac{%
(z_{1}^{2}+z_{2}^{2})b^{4}}{8}+\frac{5}{2}\right\} .  \label{e4.19**}
\end{equation}%
Replacing $t$ by $\tau $ and integrating both sides of (\ref{e4.18}) with
respect to $\tau $ over $[0,t]$, we obtain
\begin{eqnarray}
&&\lVert ^{C}\partial _{0t}^{\beta -1}u_{t}\rVert _{L^{2}(0,t;L_{\ \rho
}^{2}(\Omega ))}^{2}+\lVert ^{C}\partial _{0t}^{\gamma -1}v_{t}\rVert
_{L^{2}(0,t;L_{\ \rho }^{2}(\Omega ))}^{2}  \notag \\
&&+D_{0t}^{\beta -2}\left( \lVert u_{t}\rVert _{L_{\ \rho }^{2}(\Omega
)}^{2}+\lVert \mathcal{J}_{x}(\xi u_{t})\rVert _{L_{\ \rho }^{2}(\Omega
)}^{2}\right) +D_{0t}^{\gamma -2}\left( \lVert v_{t}\rVert _{L_{\ \rho
}^{2}(\Omega )}^{2}+\lVert \mathcal{J}_{x}(\xi v_{t})\rVert _{L_{\ \rho
}^{2}(\Omega )}^{2}\right)  \notag \\
&&+\lVert u_{x}\rVert _{L_{\ \rho }^{2}(\Omega )}^{2}+\lVert v_{x}\rVert
_{L_{\ \rho }^{2}(\Omega )}^{2}+D_{0t}^{\beta -2}\lVert u_{tx}\rVert _{L_{\
\rho }^{2}(\Omega )}^{2}+D_{0t}^{\gamma -2}\lVert u_{tx}\rVert _{L_{\ \rho
}^{2}(\Omega )}^{2}  \notag \\
&\leq &D^{\ast }\left( \int\limits_{0}^{t}\left( \lVert u_{s}\rVert _{L_{\
\rho }^{2}(\Omega )}^{2}+\lVert \mathcal{J}_{x}(\xi u_{s})\rVert _{L_{\ \rho
}^{2}(\Omega )}^{2}\right) ds+\int\limits_{0}^{t}\left( \lVert v_{s}\rVert
_{L_{\ \rho }^{2}(\Omega )}^{2}+\lVert \mathcal{J}_{x}(\xi v_{s})\rVert
_{L_{\ \rho }^{2}(\Omega )}^{2}\right) ds\right.  \notag \\
&&\left. +\int\limits_{0}^{t}\left( \lVert u_{x}\rVert _{L_{\ \rho
}^{2}(\Omega )}^{2}+\lVert v_{x}\rVert _{L_{\ \rho }^{2}(\Omega
)}^{2}\right) ds+\int\limits_{0}^{t}\left( \lVert u\rVert _{L_{\ \rho
}^{2}(\Omega )}^{2}+\lVert v\rVert _{L_{\ \rho }^{2}(\Omega )}^{2}\right)
ds\right)  \notag \\
&&+\frac{t^{2-\beta }D^{\ast }}{(2-\beta )\Gamma (2-\beta )}\left( \lVert
\varphi _{2}\rVert _{L_{\ \rho }^{2}(\Omega )}^{2}+\lVert \mathcal{J}%
_{x}(\xi \varphi _{2})\rVert _{L_{\ \rho }^{2}(\Omega )}^{2}\right) +D^{\ast
}\left( \int\limits_{0}^{t}\lVert f\rVert _{L_{\ \rho }^{2}(\Omega
)}^{2}ds+\int\limits_{0}^{t}\lVert g\rVert _{L_{\ \rho }^{2}(\Omega
)}^{2}ds\right)  \notag \\
&&+\frac{t^{2-\gamma }D^{\ast }}{(2-\gamma )\Gamma (2-\gamma )}\left( \lVert
\psi _{2}\rVert _{L_{\ \rho }^{2}(\Omega )}^{2}+\lVert \mathcal{J}_{x}(\xi
\psi _{2})\rVert _{L_{\ \rho }^{2}(\Omega )}^{2}\right) +\frac{t^{2-\beta
}D^{\ast }}{(2-\beta )\Gamma (2-\beta )}\lVert \frac{\partial \varphi _{2}}{%
\partial x}\rVert _{L_{\ \rho }^{2}(\Omega )}^{2}  \notag \\
&&+\frac{t^{2-\gamma }D^{\ast }}{(2-\gamma )\Gamma (2-\gamma )}\lVert \frac{%
\partial \psi _{2}}{\partial x}\rVert _{L_{\ \rho }^{2}(\Omega
)}^{2}+D^{\ast }\left( \lVert \frac{\partial \varphi _{1}}{\partial x}\rVert
_{L_{\ \rho }^{2}(\Omega )}^{2}+\lVert \frac{\partial \psi _{1}}{\partial x}%
\rVert _{L_{\ \rho }^{2}(\Omega )}^{2}\right) .  \label{e4.20**}
\end{eqnarray}%
Boundary integral conditions allow us to use the Poincre inequalities
\begin{equation}
\Vert u\Vert _{L_{\rho }^{2}(\Omega )}^{2}\leq \dfrac{b^{2}}{4}\left\Vert
u_{x}\right\Vert _{L_{\rho }^{2}(\Omega )}^{2},\text{ \ }\Vert v\Vert
_{L_{\rho }^{2}(\Omega )}^{2}\leq \dfrac{b^{2}}{4}\left\Vert
v_{x}\right\Vert _{L_{\rho }^{2}(\Omega )}^{2},  \label{e4.20}
\end{equation}%
to get rid of the fourth integral term on the right-hand side of (\ref%
{e4.20**}), and in the mean time, we use Poincare type inequality (\ref%
{e2.16}), we then have
\begin{eqnarray}
&&\lVert ^{C}\partial _{0t}^{\beta }u\rVert _{L^{2}(0,t;L_{\ \rho
}^{2}(\Omega ))}^{2}+\lVert ^{C}\partial _{0t}^{\gamma }v\rVert
_{L^{2}(0,t;L_{\ \rho }^{2}(\Omega ))}^{2}+\lVert ^{C}\partial _{0t}^{\beta
}u_{x}\rVert _{L^{2}(0,t;L_{\ \rho }^{2}(\Omega ))}^{2}  \notag \\
&&+\lVert ^{C}\partial _{0t}^{\gamma }v_{x}\rVert _{L^{2}(0,t;L_{\ \rho
}^{2}(\Omega ))}^{2}+D_{0t}^{\beta -2}\left( \lVert u_{t}\rVert _{L_{\ \rho
}^{2}(\Omega )}^{2}+\lVert \mathcal{J}_{x}(\xi u_{t})\rVert _{L_{\ \rho
}^{2}(\Omega )}^{2}\right)  \notag \\
&&+D_{0t}^{\gamma -2}\left( \lVert v_{t}\rVert _{L_{\ \rho }^{2}(\Omega
)}^{2}+\lVert \mathcal{J}_{x}(\xi v_{t})\rVert _{L_{\ \rho }^{2}(\Omega
)}^{2}\right) +D_{0t}^{\beta -2}\lVert u_{tx}\rVert _{L_{\ \rho }^{2}(\Omega
)}^{2}  \notag \\
&&+D_{0t}^{\gamma -2}\lVert u_{tx}\rVert _{L_{\ \rho }^{2}(\Omega
)}^{2}+\lVert u_{x}\rVert _{L_{\ \rho }^{2}(\Omega )}^{2}+\lVert v_{x}\rVert
_{L_{\ \rho }^{2}(\Omega )}^{2}  \notag \\
&\leq &D^{\ast \ast }\left( \int\limits_{0}^{t}\left( \lVert u_{s}\rVert
_{L_{\ \rho }^{2}(\Omega )}^{2}+\lVert \mathcal{J}_{x}(\xi u_{s})\rVert
_{L_{\ \rho }^{2}(\Omega )}^{2}\right) ds+\int\limits_{0}^{t}\left( \lVert
v_{s}\rVert _{L_{\ \rho }^{2}(\Omega )}^{2}+\lVert \mathcal{J}_{x}(\xi
v_{s})\rVert _{L_{\ \rho }^{2}(\Omega )}^{2}\right) ds\right.  \notag \\
&&\left. +\int\limits_{0}^{t}\left( \lVert u_{x}\rVert _{L_{\ \rho
}^{2}(\Omega )}^{2}+\lVert v_{x}\rVert _{L_{\ \rho }^{2}(\Omega
)}^{2}\right) ds\right) +D^{\ast \ast }\left( \lVert f\rVert
_{L^{2}(0,t;L_{\ \rho }^{2}(\Omega ))}^{2}+\lVert g\rVert _{L^{2}(0,t;L_{\
\rho }^{2}(\Omega ))}^{2}\right.  \notag \\
&&\left. +\lVert \varphi _{1}\rVert _{H_{\ \rho }^{1}(\Omega )}^{2}+\lVert
\psi _{1}\rVert _{H_{\ \rho }^{1}(\Omega )}^{2}+\lVert \varphi _{2}\rVert
_{H_{\ \rho }^{1}(\Omega )}^{2}+\lVert \psi _{2}\rVert _{H_{\ \rho
}^{1}(\Omega )}^{2}\right) ,  \label{e4.21***}
\end{eqnarray}%
where%
\begin{equation}
D^{\ast \ast }=D^{\ast }\max \left\{ 1,\frac{b^{4}}{2},\frac{T^{2-\beta }}{%
(2-\beta )\Gamma (2-\beta )},\frac{T^{2-\gamma }}{(2-\gamma )\Gamma
(2-\gamma )}\right\} .  \label{e4.22**}
\end{equation}%
If we leave only the last two terms on the left-hand side in inequality (\ref%
{e4.21***}), and use the Gronwall-- Bellman lemma 2.3 [55], with
\begin{eqnarray*}
R(t) &=&\int\limits_{0}^{t}\left( \lVert u_{x}\rVert _{L_{\ \rho
}^{2}(\Omega )}^{2}+\lVert v_{x}\rVert _{L_{\ \rho }^{2}(\Omega
)}^{2}\right) ds,~~R(0)=0, \\
\frac{\partial R(t)}{\partial t} &=&\lVert u_{x}\rVert _{L_{\ \rho
}^{2}(\Omega )}^{2}+\lVert v_{x}\rVert _{L_{\ \rho }^{2}(\Omega )}^{2},
\end{eqnarray*}%
we obtain%
\begin{eqnarray}
&&R(t)\leq D^{\ast \ast }e^{D^{\ast \ast }T}\left( \int\limits_{0}^{t}\left(
\lVert u_{s}\rVert _{L_{\ \rho }^{2}(\Omega )}^{2}+\lVert \mathcal{J}%
_{x}(\xi u_{s})\rVert _{L_{\ \rho }^{2}(\Omega )}^{2}\right)
ds+\int\limits_{0}^{t}\left( \lVert v_{s}\rVert _{L_{\ \rho }^{2}(\Omega
)}^{2}+\lVert \mathcal{J}_{x}(\xi v_{s})\rVert _{L_{\ \rho }^{2}(\Omega
)}^{2}\right) ds\right.  \notag \\
&&+\lVert f\rVert _{L^{2}(0,t;L_{\ \rho }^{2}(\Omega ))}^{2}+\lVert g\rVert
_{L^{2}(0,t;L_{\ \rho }^{2}(\Omega ))}^{2}+\lVert \varphi _{1}\rVert _{H_{\
\rho }^{1}(\Omega )}^{2}+\lVert \psi _{1}\rVert _{H_{\ \rho }^{1}(\Omega
)}^{2}  \label{e4.22***} \\
&&\left. +\lVert \varphi _{2}\rVert _{H_{\ \rho }^{1}(\Omega )}^{2}+\lVert
\psi _{2}\rVert _{H_{\ \rho }^{1}(\Omega )}^{2}\right) .  \notag
\end{eqnarray}%
Now by keeping only the fifth and sixth terms on the left-hand side of (\ref%
{e4.21***}), and by taking into account the inequality (\ref{e4.22***}), we
have
\begin{eqnarray}
&&D_{0t}^{\beta -2}\left( \lVert u_{t}\rVert _{L_{\ \rho }^{2}(\Omega
)}^{2}+\lVert \mathcal{J}_{x}(\xi u_{t})\rVert _{L_{\ \rho }^{2}(\Omega
)}^{2}\right) +D_{0t}^{\gamma -2}\left( \lVert v_{t}\rVert _{L_{\ \rho
}^{2}(\Omega )}^{2}+\lVert \mathcal{J}_{x}(\xi v_{t})\rVert _{L_{\ \rho
}^{2}(\Omega )}^{2}\right)  \notag \\
&\leq &\chi \left( \int\limits_{0}^{t}\left( \lVert u_{s}\rVert _{L_{\ \rho
}^{2}(\Omega )}^{2}+\lVert \mathcal{J}_{x}(\xi u_{s})\rVert _{L_{\ \rho
}^{2}(\Omega )}^{2}\right) ds+\int\limits_{0}^{t}\left( \lVert v_{s}\rVert
_{L_{\ \rho }^{2}(\Omega )}^{2}+\lVert \mathcal{J}_{x}(\xi v_{s})\rVert
_{L_{\ \rho }^{2}(\Omega )}^{2}\right) ds\right.  \notag \\
&&+\lVert f\rVert _{L^{2}(0,t;L_{\ \rho }^{2}(\Omega ))}^{2}+\lVert g\rVert
_{L^{2}(0,t;L_{\ \rho }^{2}(\Omega ))}^{2}+\lVert \varphi _{1}\rVert _{H_{\
\rho }^{1}(\Omega )}^{2}+\lVert \psi _{1}\rVert _{H_{\ \rho }^{1}(\Omega
)}^{2}  \notag \\
&&\left. +\lVert \varphi _{2}\rVert _{H_{\ \rho }^{1}(\Omega )}^{2}+\lVert
\psi _{2}\rVert _{H_{\ \rho }^{1}(\Omega )}^{2}\right) ,  \label{e4.23*}
\end{eqnarray}%
where%
\begin{equation}
\chi =D^{\ast \ast }\left( 1+D^{\ast \ast }e^{D^{\ast \ast }T}\right) .
\label{e4.24*}
\end{equation}%
By Lemma 2.1, with%
\begin{equation}
\left\{
\begin{array}{c}
\mathcal{P}_{1}(t)=\int\limits_{0}^{t}\left( \lVert u_{s}\rVert _{L_{\ \rho
}^{2}(\Omega )}^{2}+\lVert \mathcal{J}_{x}(\xi u_{s})\rVert _{L_{\ \rho
}^{2}(\Omega )}^{2}\right) ds,\text{ \ }\mathcal{P}_{1}(0)=0, \\
^{C}\partial _{0t}^{\beta -1}\mathcal{P}_{1}=D_{0t}^{\beta -2}\left( \lVert
u_{t}\rVert _{L_{\ \rho }^{2}(\Omega )}^{2}+\lVert \mathcal{J}_{x}(\xi
u_{t})\rVert _{L_{\ \rho }^{2}(\Omega )}^{2}\right) \\
\mathcal{P}_{2}(t)=\int\limits_{0}^{t}\left( \lVert v_{s}\rVert _{L_{\ \rho
}^{2}(\Omega )}^{2}+\lVert \mathcal{J}_{x}(\xi v_{s})\rVert _{L_{\ \rho
}^{2}(\Omega )}^{2}\right) ds,\text{ \ }\mathcal{P}_{2}(0)=0, \\
^{C}\partial _{0t}^{\gamma -1}\mathcal{P}_{2}=D_{0t}^{\gamma -2}\left(
\lVert v_{t}\rVert _{L_{\ \rho }^{2}(\Omega )}^{2}+\lVert \mathcal{J}%
_{x}(\xi v_{t})\rVert _{L_{\ \rho }^{2}(\Omega )}^{2}\right) ,%
\end{array}%
\right.  \label{e4.25*}
\end{equation}%
we see from (\ref{e4.23*}), that
\begin{eqnarray}
&&\mathcal{P}_{1}(t)+\mathcal{P}_{2}(t)  \notag \\
&\leq &\chi ^{\ast }\left( D^{-\beta }\lVert f\rVert _{L_{\ \rho
}^{2}(\Omega )}^{2}+D^{-\gamma }\lVert g\rVert _{L_{\ \rho }^{2}(\Omega
))}^{2}+\lVert \varphi _{1}\rVert _{H_{\ \rho }^{1}(\Omega )}^{2}+\lVert
\psi _{1}\rVert _{H_{\ \rho }^{1}(\Omega )}^{2}\right.  \notag \\
&&+\left. \lVert \varphi _{2}\rVert _{H_{\ \rho }^{1}(\Omega )}^{2}+\lVert
\psi _{2}\rVert _{H_{\ \rho }^{1}(\Omega )}^{2}\right) ,  \label{e4.26**}
\end{eqnarray}%
where
\begin{eqnarray}
\chi ^{\ast } &=&\Gamma (\beta -1)E_{\beta -1,\beta -1}(\chi t^{\beta
-1})\max \left\{ 1,\frac{T^{\beta -1}}{(\beta -1)\Gamma (\beta -1)}\right\}
\notag \\
&&+\Gamma (\gamma -1)E_{\gamma -1,\gamma -1}(\chi t^{\gamma -1})\max \left\{
1,\frac{T^{\gamma -1}}{(\gamma -1)\Gamma (\gamma -1)}\right\} .
\label{e4.27**}
\end{eqnarray}%
Owing to the inequalities%
\begin{equation}
D^{-\beta }\lVert f\rVert _{L_{\ \rho }^{2}(\Omega )}^{2}\leq \frac{t^{\beta
-1}}{\Gamma (\beta )}\int\limits_{0}^{t}\lVert f\rVert _{L_{\ \rho
}^{2}(\Omega )}^{2}ds,\text{ \ }D^{-\gamma }\lVert g\rVert _{L_{\ \rho
}^{2}(\Omega )}^{2}\leq \frac{t^{\gamma -1}}{\Gamma (\gamma )}%
\int\limits_{0}^{t}\lVert g\rVert _{L_{\ \rho }^{2}(\Omega )}^{2}ds,
\label{e4.28*}
\end{equation}%
we deduce from inequalities (\ref{e4.21***}), (\ref{e4.22***}), and (\ref%
{e4.26**}) that%
\begin{eqnarray}
&&\lVert ^{C}\partial _{0t}^{\beta }u\rVert _{L^{2}(0,t;L_{\ \rho
}^{2}(\Omega ))}^{2}+\lVert ^{C}\partial _{0t}^{\beta }u_{x}\rVert
_{L^{2}(0,t;L_{\ \rho }^{2}(\Omega ))}^{2}+\lVert ^{C}\partial _{0t}^{\gamma
}v\rVert _{L^{2}(0,t;L_{\ \rho }^{2}(\Omega ))}^{2}  \notag \\
&&+\lVert ^{C}\partial _{0t}^{\gamma }v_{x}\rVert _{L^{2}(0,t;L_{\ \rho
}^{2}(\Omega ))}^{2}+\lVert u_{x}\rVert _{L_{\ \rho }^{2}(\Omega
)}^{2}+\lVert v_{x}\rVert _{L_{\ \rho }^{2}(\Omega )}^{2}  \notag \\
&\leq &\mathcal{Y}\left( \lVert f\rVert _{L^{2}(0,t;L_{\ \rho }^{2}(\Omega
))}^{2}+\lVert g\rVert _{L^{2}(0,t;L_{\ \rho }^{2}(\Omega ))}^{2}+\lVert
\varphi _{1}\rVert _{H_{\ \rho }^{1}(\Omega )}^{2}+\lVert \psi _{1}\rVert
_{H_{\ \rho }^{1}(\Omega )}^{2}\right.  \notag \\
&&+\left. \lVert \varphi _{2}\rVert _{H_{\ \rho }^{1}(\Omega )}^{2}+\lVert
\psi _{2}\rVert _{H_{\ \rho }^{1}(\Omega )}^{2}\right) ,  \label{e4.29**}
\end{eqnarray}%
where%
\begin{equation}
\mathcal{Y=}\chi ^{\ast }\chi \max \left( \frac{T^{\beta -1}}{\Gamma (\beta )%
},\frac{T^{\gamma -1}}{\Gamma (\gamma )}\right) .  \label{e4.29***}
\end{equation}%
By virtue of poincare inequalities (\ref{e4.20}), and equivalence of norms,
the inequality (\ref{e4.29**}) takes the form
\begin{eqnarray}
&&\lVert ^{C}\partial _{0t}^{\beta }u\rVert _{L^{2}(0,t;H_{\ \rho
}^{1}(\Omega ))}^{2}+\lVert ^{C}\partial _{0t}^{\gamma }v\rVert
_{L^{2}(0,t;H_{\ \rho }^{1}(\Omega ))}^{2}+\lVert u\rVert _{H_{\ \rho
}^{1}(\Omega )}^{2}+\lVert v\rVert _{H_{\ \rho }^{1}(\Omega )}^{2}  \notag \\
&\leq &\mathcal{Y}^{\ast }\left( \lVert f\rVert _{L^{2}(0,t;L_{\ \rho
}^{2}(\Omega ))}^{2}+\lVert g\rVert _{L^{2}(0,t;L_{\ \rho }^{2}(\Omega
))}^{2}+\lVert \varphi _{1}\rVert _{H_{\ \rho }^{1}(\Omega )}^{2}+\lVert
\psi _{1}\rVert _{H_{\ \rho }^{1}(\Omega )}^{2}\right.  \notag \\
&&+\left. \lVert \varphi _{2}\rVert _{H_{\ \rho }^{1}(\Omega )}^{2}+\lVert
\psi _{2}\rVert _{H_{\ \rho }^{1}(\Omega )}^{2}\right) ,  \label{e4.30}
\end{eqnarray}%
where%
\begin{equation}
\mathcal{Y}^{\ast }=\frac{\mathcal{Y}}{\min \left( 1,\frac{b^{2}}{4}\right) }%
.  \label{e4.30**}
\end{equation}%
Now by adding the quantity $\lVert u\rVert _{L^{2}(0,t;H_{\ \rho
}^{1}(\Omega ))}^{2}+\lVert v\rVert _{L^{2}(0,t;H_{\ \rho }^{1}(\Omega
))}^{2}$ to both sides of (\ref{e4.30}), we have%
\begin{eqnarray}
&&\lVert ^{C}\partial _{0t}^{\beta }u\rVert _{L^{2}(0,t;H_{\ \rho
}^{1}(\Omega ))}^{2}+\lVert u\rVert _{L^{2}(0,t;H_{\ \rho }^{1}(\Omega
))}^{2}+\lVert u\rVert _{H_{\ \rho }^{1}(\Omega )}^{2}  \notag \\
&&+\lVert ^{C}\partial _{0t}^{\gamma }v\rVert _{L^{2}(0,t;H_{\ \rho
}^{1}(\Omega ))}^{2}+\lVert v\rVert _{L^{2}(0,t;H_{\ \rho }^{1}(\Omega
))}^{2}+\lVert v\rVert _{H_{\ \rho }^{1}(\Omega )}^{2}  \notag \\
&\leq &\mathcal{Y}^{\ast \ast }\left( \lVert f\rVert _{L^{2}(0,t;L_{\ \rho
}^{2}(\Omega ))}^{2}+\lVert g\rVert _{L^{2}(0,t;L_{\ \rho }^{2}(\Omega
))}^{2}+\lVert \varphi _{1}\rVert _{H_{\ \rho }^{1}(\Omega )}^{2}+\lVert
\psi _{1}\rVert _{H_{\ \rho }^{1}(\Omega )}^{2}\right.  \notag \\
&&+\left. \lVert \varphi _{2}\rVert _{H_{\ \rho }^{1}(\Omega )}^{2}+\lVert
\psi _{2}\rVert _{H_{\ \rho }^{1}(\Omega )}^{2}+\lVert u\rVert
_{L^{2}(0,t;H_{\ \rho }^{1}(\Omega ))}^{2}+\lVert v\rVert _{L^{2}(0,t;H_{\
\rho }^{1}(\Omega ))}^{2}\right)  \label{e4.31}
\end{eqnarray}%
\begin{equation}
\mathcal{Y}^{\ast \ast }=\max \left( 1,\mathcal{Y}^{\ast }\right) .
\label{e4.32*}
\end{equation}%
Application of Gronwall's Lemma to (\ref{e4.31}) gives the inequality%
\begin{eqnarray}
&&\lVert u\rVert _{\mathcal{W}^{\beta }(Q_{t})}^{2}+\lVert v\rVert _{%
\mathcal{W}^{\gamma }(Q_{t})}^{2}+\lVert u\rVert _{H_{\ \rho }^{1}(\Omega
)}^{2}+\lVert v\rVert _{H_{\ \rho }^{1}(\Omega )}^{2}  \notag \\
&\leq &\mathcal{Y}^{\ast \ast }e^{T\mathcal{Y}^{\ast \ast }}\left( \lVert
f\rVert _{L^{2}(0,T;L_{\ \rho }^{2}(\Omega ))}^{2}+\lVert g\rVert
_{L^{2}(0,T;L_{\ \rho }^{2}(\Omega ))}^{2}+\lVert \varphi _{1}\rVert _{H_{\
\rho }^{1}(\Omega )}^{2}+\lVert \psi _{1}\rVert _{H_{\ \rho }^{1}(\Omega
)}^{2}\right.  \notag \\
&&+\left. \lVert \varphi _{2}\rVert _{H_{\ \rho }^{1}(\Omega )}^{2}+\lVert
\psi _{2}\rVert _{H_{\ \rho }^{1}(\Omega )}^{2}\right) .  \label{e4.33*}
\end{eqnarray}

The independence of the right-hand side on $t$ in (\ref{e4.33*}), gives%
\begin{eqnarray}
&&\lVert u\rVert _{\mathcal{W}^{\beta }(Q_{T})}^{2}+\lVert u\rVert _{%
\mathcal{W}^{\gamma }(Q_{T})}^{2}+\lVert u\rVert _{C(0,T.,H_{\ \rho
}^{1}(\Omega )}^{2}+\lVert v\rVert _{C(0,T.,H_{\ \rho }^{1}(\Omega )}^{2}
\notag \\
&\leq &\mathcal{M}\left( \lVert f\rVert _{L^{2}(0,T;L_{\ \rho }^{2}(\Omega
))}^{2}+\lVert g\rVert _{L^{2}(0,T;L_{\ \rho }^{2}(\Omega ))}^{2}+\lVert
\varphi _{1}\rVert _{H_{\ \rho }^{1}(\Omega )}^{2}+\lVert \psi _{1}\rVert
_{H_{\ \rho }^{1}(\Omega )}^{2}\right.  \notag \\
&&+\left. \lVert \varphi _{2}\rVert _{H_{\ \rho }^{1}(\Omega )}^{2}+\lVert
\psi _{2}\rVert _{H_{\ \rho }^{1}(\Omega )}^{2}\right) ,  \label{e4.34*}
\end{eqnarray}%
where $\mathcal{M=Y}^{\ast \ast }e^{T\mathcal{Y}^{\ast \ast }}.$

It can be proved in a standard way that the operator $\mathcal{X}%
:B\rightarrow H$ is closable. Let $\overline{\mathcal{X}}$ be its closure.

\begin{prop}
The operator $\mathcal{X}:B\rightarrow H$ has a closure.
\end{prop}

\textbf{Proof:} The proof can be established in a similar way as in [57].

These are some consequences of Theorem 4.1.

\begin{coro}
\label{coro4.1} There exists a positive constant $C$ such that
\begin{equation}
\Vert W\Vert _{B}\leq C\Vert \overline{\mathcal{X}}W\Vert _{H},\quad \forall
W\in D(\overline{\mathcal{X}}),  \label{e4.28}
\end{equation}%
where: $C=\sqrt{C_{7}}$. \newline
The inequality (\ref{e4.28}) means that inequality (\ref{e4.1}) can be
extended to strong solutions after passing to limit.
\end{coro}

We can deduce from inequality (\ref{e4.28}) that a strong solution of the
system (\ref{e3.1})-(\ref{e3.3}) is unique and depends continuously on $%
\mathcal{F}=\left( \mathcal{F}_{1},\mathcal{F}_{2}\right) \in H,$ where $%
\mathcal{F}_{1}=\left\{ f,\varphi _{1},\varphi _{2}\right\} $ and $\mathcal{F%
}_{2}=\left\{ g,\psi _{1},\psi _{2}\right\} $, and that the image $R(%
\overline{\mathcal{X}})$ of $\overline{\mathcal{X}}$ is closed in $H$ and $R(%
\overline{\mathcal{X}})=\overline{R(\mathcal{X})}$. So in order to prove
that the system (\ref{e3.1})-(\ref{e3.3}) has a strong solution for
arbitrary $\left( \mathcal{F}_{1},\mathcal{F}_{2}\right) \in H,$ it is
sufficient to prove that the range of $\mathcal{X}$ is dense in $H$, that is
$\overline{R(\mathcal{X})}=H$.

\section{Existence of the solution of the linear system}

\begin{prop}
\label{prop5.1} If for some function: $Y^{\ast }(x,t)=(y_{1}^{\ast
}(x,t),y_{2}^{\ast }(x,t))\in ($\ $L^{2}(0,T;L_{\rho }^{2}(\Omega )))^{2}$,
and for all $W(x,t)=(u(x,t),v(x,t))\in D_{0}(\mathcal{X})=\left\{ W/W\in D(%
\mathcal{X}):\ell _{1}u=0,\ell _{2}u=0,\ell _{3}v=0,\ell _{4}v=0\right\} $,
we have
\begin{equation}
\left( \mathcal{L}W,Y^{\ast }\right) _{\ L^{2}(0,T;L_{\rho }^{2}(\Omega
))}=\left( \mathcal{L}_{1}u,y_{1}^{\ast }\right) _{L^{2}(0,T;L_{\rho
}^{2}(\Omega ))}+\left( \mathcal{L}_{2}v,y_{2}^{\ast }\right)
_{L^{2}(0,T;L_{\rho }^{2}(\Omega ))}=0,  \label{e5.1}
\end{equation}%
then $Y^{\ast }$ vanishes a.e in the domain $Q$.
\end{prop}

\textbf{Proof:} We first set
\begin{eqnarray}
&&W=(u,v)=(\mathcal{J}_{t}^{2}(p_{1}),\mathcal{J}_{t}^{2}(p_{2})),
\label{e5.2} \\
&&Y^{\ast }(x,t)=(y_{1}^{\ast }(x,t),y_{2}^{\ast }(x,t))=(\mathcal{J}%
_{t}(p_{1})-\mathcal{J}_{x}^{2}\left( \xi \mathcal{J}_{t}(p_{1})\right) ,%
\mathcal{J}_{t}(p_{2})-\mathcal{J}_{x}^{2}\left( \xi \mathcal{J}%
_{t}(p_{2}\right) ),  \label{e5.3}
\end{eqnarray}%
where
\begin{eqnarray*}
\mathcal{J}_{t}(p_{i}) &=&\int\limits_{0}^{t}p_{i}(x,s)ds,\text{ \ }\mathcal{%
J}_{t}^{2}(p_{i})=\int\limits_{0}^{t}\int\limits_{0}^{s}p_{i}(x,z)dzds, \\
\mathcal{J}_{x}^{2}\left( \xi \mathcal{J}_{t}(p_{i})\right)
&=&\int\limits_{0}^{x}\int\limits_{0}^{\xi }\int\limits_{0}^{t}\eta
~p_{i}(\eta ,s)dsd\eta d\xi ,~~~i=1,2.
\end{eqnarray*}%
We suppose that the functions $~~p_{i}(x,t)$ satisfy conditions (\ref{e3.3})
and such that
\begin{equation*}
p_{i},\text{ }p_{ix},\text{ }\mathcal{J}_{t}(p_{i}),\text{ }\mathcal{J}%
_{t}^{2}(p_{i}),\text{ }x\mathcal{J}_{t}^{2}(p_{ix}),\text{ }\mathcal{J}%
_{x}^{2}\left( \xi \mathcal{J}_{t}(p_{i})\right) ,\text{ }^{C}\partial
_{0t}^{\beta }p_{i},\text{ }^{C}\partial _{0t}^{\gamma }p_{i}\in L^{2}(Q),%
\text{ \ }i=1,2.
\end{equation*}%
Now by replacing (\ref{e5.2}) and (\ref{e5.3}) in the relation (\ref{e5.1}),
we obtain
\begin{eqnarray}
&&\left( ^{C}\partial _{0t}^{\beta }\left( \mathcal{J}_{t}^{2}(p_{1})\right)
,\mathcal{J}_{t}(p_{1})\right) _{{L_{\ \rho }^{2}\left( \Omega \right) }%
}-\left( \left( x\left( \mathcal{J}_{t}^{2}(p_{1x})\right) \right) _{x},%
\mathcal{J}_{t}(p_{1})\right) _{{L^{2}\left( \Omega \right) }}-\left( \left(
x\left( \mathcal{J}_{t}^{2}(p_{1x})\right) \right) _{xt},\mathcal{J}%
_{t}(p_{1})\right) _{{L^{2}\left( \Omega \right) }}  \notag \\
&&+\left( \mathcal{J}_{t}^{2}(p_{2}),\mathcal{J}_{t}(p_{1})\right) _{L_{\
\rho }^{2}\left( \Omega \right) }-(^{C}\partial _{0t}^{\beta }\left(
\mathcal{J}_{t}^{2}(p_{1})\right) ,\mathcal{J}_{x}^{2}\left( \xi \mathcal{J}%
_{t}(p_{1})\right) )_{{L_{\ \rho }^{2}\left( \Omega \right) }}+\left( \left(
x\left( \mathcal{J}_{t}^{2}(p_{1x})\right) \right) _{x},\mathcal{J}%
_{x}^{2}\left( \xi \mathcal{J}_{t}(p_{1})\right) \right) _{{L^{2}\left(
\Omega \right) }}  \notag \\
&&+\left( \left( x\left( \mathcal{J}_{t}^{2}(p_{1x})\right) \right) _{xt},%
\mathcal{J}_{x}^{2}\left( \xi \mathcal{J}_{t}(p_{1})\right) \right) _{{%
L^{2}\left( \Omega \right) }}-\left( \mathcal{J}_{t}^{2}(p_{2}),\mathcal{J}%
_{x}^{2}\left( \xi \mathcal{J}_{t}(p_{1})\right) \right) _{L_{\ \rho
}^{2}\left( \Omega \right) }  \notag \\
&&+\left( ^{C}\partial _{0t}^{\gamma }\left( \mathcal{J}_{t}^{2}(p_{2})%
\right) ,\mathcal{J}_{t}(p_{2})\right) _{{L_{\ \rho }^{2}\left( \Omega
\right) }}-\left( \left( x\left( \mathcal{J}_{t}^{2}(p_{2x})\right) \right)
_{x},\mathcal{J}_{t}(p_{2})\right) _{{L^{2}\left( \Omega \right) }}-\left(
\left( x\left( \mathcal{J}_{t}^{2}(p_{2x})\right) \right) _{xt},\mathcal{J}%
_{t}(p_{2})\right) _{{L^{2}\left( \Omega \right) }}  \notag \\
&&+\left( \mathcal{J}_{t}^{2}(p_{1}),\mathcal{J}_{t}(p_{2})\right) _{L_{\
\rho }^{2}\left( \Omega \right) }-\left( ^{C}\partial _{0t}^{\gamma }\left(
\mathcal{J}_{t}^{2}(p_{2})\right) ,\mathcal{J}_{x}^{2}\left( \xi \mathcal{J}%
_{t}(p_{2})\right) \right) _{{L_{\ \rho }^{2}\left( \Omega \right) }}+\left(
\left( x\left( \mathcal{J}_{t}^{2}(p_{2x})\right) \right) _{x},\mathcal{J}%
_{x}^{2}\left( \xi \mathcal{J}_{t}(p_{2})\right) \right) _{{L^{2}\left(
\Omega \right) }}  \notag \\
&&+\left( \left( x\left( \mathcal{J}_{t}^{2}(p_{2x})\right) \right) _{xt},%
\mathcal{J}_{x}^{2}\left( \xi \mathcal{J}_{t}(p_{2})\right) \right) _{{%
L^{2}\left( \Omega \right) }}-\left( \mathcal{J}_{t}^{2}(p_{1}),\mathcal{J}%
_{x}^{2}\left( \xi \mathcal{J}_{t}(p_{2})\right) \right) _{L_{\ \rho
}^{2}\left( \Omega \right) }=0.  \label{e5.4}
\end{eqnarray}%
Since%
\begin{equation*}
\lVert \mathcal{J}_{t}^{2}(p_{i})\rVert _{L_{\ \rho }^{2}(\Omega )}^{2}\leq
\dfrac{T^{2}}{2}\lVert \mathcal{J}_{t}(p_{i})\rVert _{L_{\ \rho }^{2}(\Omega
)}^{2},~i=1,2,
\end{equation*}%
then, using conditions (\ref{e3.3}), and computation of each term of (\ref%
{e5.4}), gives%
\begin{eqnarray}
\left( ^{C}\partial _{0t}^{\beta }\left( \mathcal{J}_{t}^{2}(p_{1})\right) ,%
\mathcal{J}_{t}(p_{1})\right) _{{L_{\ \rho }^{2}\left( \Omega \right) }}
&=&\left( ^{C}\partial _{0t}^{\beta -1}\left( \mathcal{J}_{t}(p_{1})\right) ,%
\mathcal{J}_{t}(p_{1})\right) _{{L_{\ \rho }^{2}\left( \Omega \right) }}
\notag \\
&\geq &\dfrac{1}{2}~~^{C}\partial _{0t}^{\beta -1}\lVert \mathcal{J}%
_{t}(p_{1})\rVert _{L_{\ \rho }^{2}(\Omega )}^{2},  \label{e5.5}
\end{eqnarray}%
\begin{equation}
-\left( \left( x\left( \mathcal{J}_{t}^{2}(p_{1x})\right) \right) _{x},%
\mathcal{J}_{t}(p_{1})\right) _{{L^{2}\left( \Omega \right) }}=\dfrac{1}{2}%
\dfrac{\partial }{\partial t}\lVert \mathcal{J}_{t}^{2}(p_{1x})\rVert _{L_{\
\rho }^{2}(\Omega )}^{2},  \label{e5.6}
\end{equation}%
\begin{equation}
-\left( \left( x\left( \mathcal{J}_{t}^{2}(p_{1x})\right) \right) _{xt},%
\mathcal{J}_{t}(p_{1})\right) _{{L^{2}\left( \Omega \right) }}=\lVert
\mathcal{J}_{t}(p_{1x})\rVert _{L_{\ \rho }^{2}(\Omega )}^{2},  \label{e5.7}
\end{equation}%
\begin{eqnarray}
-\left( \mathcal{J}_{t}^{2}(p_{2}),\mathcal{J}_{t}(p_{1})\right) _{L_{\ \rho
}^{2}\left( \Omega \right) } &\leq &\dfrac{1}{2}\lVert \mathcal{J}%
_{t}^{2}(p_{2})\rVert _{L_{\ \rho }^{2}(\Omega )}^{2}+\dfrac{1}{2}\lVert
\mathcal{J}_{t}(p_{1})\rVert _{L_{\ \rho }^{2}(\Omega )}^{2}  \notag \\
&\leq &\dfrac{T^{2}}{4}\lVert \mathcal{J}_{t}(p_{2})\rVert _{L_{\ \rho
}^{2}(\Omega )}^{2}+\dfrac{1}{2}\lVert \mathcal{J}_{t}(p_{1})\rVert _{L_{\
\rho }^{2}(\Omega )}^{2},  \label{e5.8}
\end{eqnarray}%
\begin{eqnarray}
-\left( ^{C}\partial _{0t}^{\beta }\left( \mathcal{J}_{t}^{2}(p_{1})\right) ,%
\mathcal{J}_{x}^{2}\left( \xi \mathcal{J}_{t}(p_{1})\right) \right) _{{L_{\
\rho }^{2}\left( \Omega \right) }} &=&-\left( ^{C}\partial _{0t}^{\beta
-1}\left( \mathcal{J}_{t}(p_{1})\right) ,\mathcal{J}_{x}^{2}\left( \xi
\mathcal{J}_{t}(p_{1})\right) \right) _{{L_{\ \rho }^{2}\left( \Omega
\right) }}  \notag \\
&=&\left( ^{C}\partial _{0t}^{\beta -1}\left( \mathcal{J}_{x}\left( \xi
\mathcal{J}_{t}(p_{1})\right) \right) ,\mathcal{J}_{x}\left( \xi \mathcal{J}%
_{t}(p_{1})\right) \right) _{{L^{2}\left( \Omega \right) }}  \notag \\
&\geq &\dfrac{1}{2b}~~^{C}\partial _{0t}^{\beta -1}\lVert \mathcal{J}%
_{x}\left( \xi \mathcal{J}_{t}(p_{1})\right) \rVert _{L_{\ \rho }^{2}(\Omega
)}^{2},  \label{e5.9}
\end{eqnarray}%
\begin{eqnarray}
\left( \left( x\left( \mathcal{J}_{t}^{2}(p_{1x})\right) \right) _{x},%
\mathcal{J}_{x}^{2}\left( \xi \mathcal{J}_{t}(p_{1})\right) \right) _{{%
L^{2}\left( \Omega \right) }} &=&-\left( \mathcal{J}_{t}^{2}(p_{1x}),%
\mathcal{J}_{x}\left( \xi \mathcal{J}_{t}(p_{1})\right) \right) _{{L_{\ \rho
}^{2}\left( \Omega \right) }}  \notag \\
&\leq &\dfrac{1}{T^{2}}\lVert \mathcal{J}_{t}^{2}(p_{1x})\rVert _{L_{\ \rho
}^{2}(\Omega )}^{2}+\dfrac{T^{2}}{4}\lVert \mathcal{J}_{x}\left( \xi
\mathcal{J}_{t}(p_{1})\right) \rVert _{L_{\ \rho }^{2}(\Omega )}^{2}  \notag
\\
&\leq &\dfrac{1}{2}\lVert \mathcal{J}_{t}(p_{1x})\rVert _{L_{\ \rho
}^{2}(\Omega )}^{2}+\dfrac{T^{2}}{4}\lVert \mathcal{J}_{x}\left( \xi
\mathcal{J}_{t}(p_{1})\right) \rVert _{L_{\ \rho }^{2}(\Omega )}^{2},
\label{e5.10}
\end{eqnarray}%
\begin{eqnarray}
\left( \left( x\left( \mathcal{J}_{t}^{2}(p_{1x})\right) \right) _{xt},%
\mathcal{J}_{x}^{2}\left( \xi \mathcal{J}_{t}(p_{1})\right) \right) _{{%
L^{2}\left( \Omega \right) }} &=&-\left( \mathcal{J}_{t}(p_{1x}),\mathcal{J}%
_{x}\left( \xi \mathcal{J}_{t}(p_{1})\right) \right) _{{L_{\ \rho
}^{2}\left( \Omega \right) }}  \notag \\
&\leq &\dfrac{1}{2}\lVert \mathcal{J}_{t}(p_{1x})\rVert _{L_{\ \rho
}^{2}(\Omega )}^{2}+\dfrac{1}{2}\lVert \mathcal{J}_{x}\left( \xi \mathcal{J}%
_{t}(p_{1})\right) \rVert _{L_{\ \rho }^{2}(\Omega )}^{2},  \label{e5.11}
\end{eqnarray}%
\begin{eqnarray}
\left( \mathcal{J}_{t}^{2}(p_{2}),\mathcal{J}_{x}^{2}\left( \xi \mathcal{J}%
_{t}(p_{1})\right) \right) _{L_{\ \rho }^{2}\left( \Omega \right) } &\leq &%
\dfrac{1}{2}\lVert \mathcal{J}_{t}^{2}(p_{2})\rVert _{L_{\ \rho }^{2}(\Omega
)}^{2}+\dfrac{1}{2}\lVert \mathcal{J}_{x}^{2}\left( \xi \mathcal{J}%
_{t}(p_{1})\right) \rVert _{L_{\ \rho }^{2}(\Omega )}^{2}  \notag \\
&\leq &\dfrac{T^{2}}{4}\lVert \mathcal{J}_{t}(p_{2})\rVert _{L_{\ \rho
}^{2}(\Omega )}^{2}+\dfrac{b^{6}}{8}\lVert \mathcal{J}_{t}(p_{1})\rVert
_{L_{\ \rho }^{2}(\Omega )}^{2}.  \label{e5.12}
\end{eqnarray}%
Combination of (\ref{e5.5})--(\ref{e5.12}) and (\ref{e5.4}), yields
\begin{eqnarray}
&&~^{C}\partial _{0t}^{\beta -1}\left( \lVert \mathcal{J}_{t}(p_{1})\rVert
_{L_{\ \rho }^{2}(\Omega )}^{2}+\lVert \mathcal{J}_{x}\left( \xi \mathcal{J}%
_{t}(p_{1})\right) \rVert _{L_{\ \rho }^{2}(\Omega )}^{2}\right) +\dfrac{%
\partial }{\partial t}\lVert \mathcal{J}_{t}^{2}(p_{1x})\rVert _{L_{\ \rho
}^{2}(\Omega )}^{2}  \notag \\
&&~^{C}\partial _{0t}^{\gamma -1}\left( \lVert \mathcal{J}_{t}(p_{2})\rVert
_{L_{\ \rho }^{2}(\Omega )}^{2}+\lVert \mathcal{J}_{x}\left( \xi \mathcal{J}%
_{t}(p_{2})\right) \rVert _{L_{\ \rho }^{2}(\Omega )}^{2}\right) +\dfrac{%
\partial }{\partial t}\lVert \mathcal{J}_{t}^{2}(p_{2x})\rVert _{L_{\ \rho
}^{2}(\Omega )}^{2}  \notag \\
&\leq &M_{1}\left( \lVert \mathcal{J}_{t}(p_{1})\rVert _{L_{\ \rho
}^{2}(\Omega )}^{2}+\lVert \mathcal{J}_{x}\left( \xi \mathcal{J}%
_{t}(p_{1})\right) \rVert _{L_{\ \rho }^{2}(\Omega )}^{2}\right.  \notag \\
&&+\left. \lVert \mathcal{J}_{t}(p_{2})\rVert _{L_{\ \rho }^{2}(\Omega
)}^{2}+\lVert \mathcal{J}_{x}\left( \xi \mathcal{J}_{t}(p_{2})\right) \rVert
_{L_{\ \rho }^{2}(\Omega )}^{2}\right) ,  \label{e5.13*}
\end{eqnarray}%
where
\begin{equation*}
M_{1}=\dfrac{\max \left\{ T^{2},1+\frac{b^{6}}{4},1+\frac{T^{2}}{2}\right\}
}{\min \left\{ 1,\frac{1}{b}\right\} }.
\end{equation*}%
After integration, we entail from (\ref{e5.13*}) that
\begin{eqnarray}
&&D_{0t}^{\beta -2}\lVert \mathcal{J}_{t}(p_{1})\rVert _{L_{\ \rho
}^{2}(\Omega )}^{2}+D_{0t}^{\beta -2}\lVert \mathcal{J}_{x}\left( \xi
\mathcal{J}_{t}(p_{1})\right) \rVert _{L_{\ \rho }^{2}(\Omega )}^{2}+\lVert
\mathcal{J}_{t}^{2}(p_{1x})\rVert _{L_{\ \rho }^{2}(\Omega )}^{2}  \notag \\
&&D_{0t}^{\gamma -2}\lVert \mathcal{J}_{t}(p_{2})\rVert _{L_{\ \rho
}^{2}(\Omega )}^{2}+D_{0t}^{\gamma -2}\lVert \mathcal{J}_{x}\left( \xi
\mathcal{J}_{t}(p_{2})\right) \rVert _{L_{\ \rho }^{2}(\Omega )}^{2}+\lVert
\mathcal{J}_{t}^{2}(p_{2x})\rVert _{L_{\ \rho }^{2}(\Omega )}^{2}  \notag \\
&\leq &M_{1}\left[ \int\limits_{0}^{t}\left( \lVert \mathcal{J}_{\tau
}(p_{1})\rVert _{L_{\ \rho }^{2}(\Omega )}^{2}+\lVert \mathcal{J}_{x}\left(
\xi \mathcal{J}_{\tau }(p_{1})\right) \rVert _{L_{\ \rho }^{2}(\Omega
)}^{2}\right) d\tau \right.  \notag \\
&&+\left. \int\limits_{0}^{t}\left( \lVert \mathcal{J}_{\tau }(p_{2})\rVert
_{L_{\ \rho }^{2}(\Omega )}^{2}+\lVert \mathcal{J}_{x}\left( \xi \mathcal{J}%
_{\tau }(p_{2})\right) \rVert _{L_{\ \rho }^{2}(\Omega )}^{2}\right) d\tau %
\right] .  \label{e5.14}
\end{eqnarray}%
If we drop the last four terms on the left-hand side of (\ref{e5.14}), apply
Lemma 2.1, and use inequality (\ref{e2.14}) , we have
\begin{eqnarray}
&&\int\limits_{0}^{t}\left( \lVert \mathcal{J}_{\tau }(p_{1})\rVert _{L_{\
\rho }^{2}(\Omega )}^{2}+\lVert \mathcal{J}_{x}\left( \xi \mathcal{J}_{\tau
}(p_{1})\right) \rVert _{L_{\ \rho }^{2}(\Omega )}^{2}\right) d\tau  \notag
\\
&\leq &M_{1}\Gamma (\beta -1)E_{\beta -1,\beta -1}(M_{1}T^{\beta
-1})D_{0t}^{-\beta }\left( \lVert \mathcal{J}_{\tau }(p_{2})\rVert _{L_{\
\rho }^{2}(\Omega )}^{2}+\lVert \mathcal{J}_{x}\left( \xi \mathcal{J}_{\tau
}(p_{2})\right) \rVert _{L_{\ \rho }^{2}(\Omega )}^{2}\right) .
\label{e5.15*}
\end{eqnarray}%
Application of inequality (\ref{e2.14}), reduces (\ref{e5.15*}) to%
\begin{eqnarray}
&&\int\limits_{0}^{t}\left( \lVert \mathcal{J}_{\tau }(p_{1})\rVert _{L_{\
\rho }^{2}(\Omega )}^{2}+\lVert \mathcal{J}_{x}\left( \xi \mathcal{J}_{\tau
}(p_{1})\right) \rVert _{L_{\ \rho }^{2}(\Omega )}^{2}\right) d\tau  \notag
\\
&\leq &M_{2}\left( \int\limits_{0}^{t}\left( \lVert \mathcal{J}_{\tau
}(p_{2})\rVert _{L_{\ \rho }^{2}(\Omega )}^{2}+\lVert \mathcal{J}_{x}\left(
\xi \mathcal{J}_{\tau }(p_{2})\right) \rVert _{L_{\ \rho }^{2}(\Omega
)}^{2}\right) d\tau \right) ,  \label{e5.15**}
\end{eqnarray}%
where
\begin{equation}
M_{2}=M_{1}\Gamma (\beta -1)E_{\beta -1,\beta -1}(M_{1}T^{\beta -1})\dfrac{%
T^{\beta -1}}{\Gamma (\beta )}.  \label{e5.15***}
\end{equation}%
We infer from inequalities (\ref{e5.15**}) and (\ref{e5.14}) that
\begin{eqnarray}
&&D_{0t}^{\gamma -2}\lVert \mathcal{J}_{t}(p_{2})\rVert _{L_{\ \rho
}^{2}(\Omega )}^{2}+D_{0t}^{\gamma -2}\lVert \mathcal{J}_{x}\left( \xi
\mathcal{J}_{t}(p_{2})\right) \rVert _{L_{\ \rho }^{2}(\Omega )}^{2}+\lVert
\mathcal{J}_{t}^{2}(p_{2x})\rVert _{L_{\ \rho }^{2}(\Omega )}^{2}  \notag \\
&&D_{0t}^{\beta -2}\lVert \mathcal{J}_{t}(p_{1})\rVert _{L_{\ \rho
}^{2}(\Omega )}^{2}+D_{0t}^{\beta -2}\lVert \mathcal{J}_{x}\left( \xi
\mathcal{J}_{t}(p_{1})\right) \rVert _{L_{\ \rho }^{2}(\Omega )}^{2}+\lVert
\mathcal{J}_{t}^{2}(p_{1x})\rVert _{L_{\ \rho }^{2}(\Omega )}^{2}  \notag \\
&\leq &M_{3}\left[ \int\limits_{0}^{t}\left( \lVert \mathcal{J}_{\tau
}(p_{2})\rVert _{L_{\ \rho }^{2}(\Omega )}^{2}+\lVert \mathcal{J}_{x}\left(
\xi \mathcal{J}_{\tau }(p_{2})\right) \rVert _{L_{\ \rho }^{2}(\Omega
)}^{2}\right) d\tau \right] ,  \label{e5.16}
\end{eqnarray}%
where
\begin{equation*}
M_{3}=M_{1}(1+M_{2}).
\end{equation*}%
If we now discard the last four terms in the left-hand side of (\ref{e5.16}%
), and apply Lemma 2.1, we get
\begin{equation*}
\int\limits_{0}^{t}\left( \lVert \mathcal{J}_{\tau }(p_{2})\rVert _{L_{\
\rho }^{2}(\Omega )}^{2}+\lVert \mathcal{J}_{x}\left( \xi \mathcal{J}_{\tau
}(p_{2})\right) \rVert _{L_{\ \rho }^{2}(\Omega )}^{2}\right) d\tau \leq
M_{4}\left( D_{0t}^{-\gamma }(0)\right) =0,
\end{equation*}%
with $M_{4}=\Gamma (\gamma -1)E_{\gamma -1,\gamma -1}(M_{3}T^{\gamma -1})$.%
\newline
Hence, we deduce that $Y^{\ast }(x,t)=(y_{1}^{\ast },y_{2}^{\ast })=(0,0)$
almost everywhere in the domain $Q$.

\begin{thm}
\label{thm5.1} For any $(f,g)\in \left( L_{\rho }^{2}(Q)\right) ^{2}$ and
any $(\varphi _{1},\psi _{1}),(\varphi _{2},\psi _{2})\in \left( H_{\rho
}^{1}(\Omega )\right) ^{2},$ there exists a unique strong solution $W=%
\overline{\mathcal{X}}^{-1}\mathcal{F}=\overline{\mathcal{X}^{-1}}\mathcal{F}
$ of the system (\ref{e3.1})-(\ref{e3.3}), where $\mathcal{F}=\left(
\mathcal{F}_{1},\mathcal{F}_{2}\right) \in H,~~~\mathcal{F}_{1}=\left\{
f,\varphi _{1},\varphi _{2}\right\} ,~~\mathcal{F}_{2}=\left\{ g,\psi
_{1},\psi _{2}\right\} ,~~W=(u,v)$ and
\begin{equation*}
\Vert W\Vert _{B}\leq C\Vert \mathcal{X}W\Vert _{H},
\end{equation*}%
for a positive constant $C,$ independent of $W$.
\end{thm}

\textbf{Proof:} We show the validity of $\overline{R(\mathcal{X})}=H.$ Since
$H$ is a Hilbert space, the equality $\overline{R(\mathcal{X})}=H$ holds if
\begin{eqnarray}
(LW,Y)_{H} &=&\left( \left\{ L_{1}(u,v),L_{2}(u,v)\right\} ,\left\{ \mathcal{%
Y}_{1},\mathcal{Y}_{2}\right\} \right) _{H}  \notag \\
&=&\left( \left\{ \left( \mathcal{L}_{1}(u,v),\ell _{1}u,\ell _{2}u\right)
,\left( \mathcal{L}_{2}(u,v),\ell _{3}v,\ell _{4}v\right) \right\} ,\left\{
\left( y_{1},y_{2},y_{3}\right) ,\left( y_{4},y_{5},y_{6}\right) \right\}
\right) _{H}  \notag \\
&=&\left( \mathcal{L}_{1}(u,v),y_{1}\right) _{L^{2}(0,T;L_{\rho }^{2}(\Omega
))}+\left( \ell _{1}u,y_{2}\right) _{H_{\rho }^{1}(\Omega )}+\left( \ell
_{2}u,y_{3}\right) _{H_{\rho }^{1}(\Omega )}  \notag \\
&&+\left( \mathcal{L}_{2}(u,v),y_{4}\right) _{L^{2}(0,T;L_{\rho }^{2}(\Omega
))}+\left( \ell _{3}v,y_{5}\right) _{H_{\rho }^{1}(\Omega )}+\left( \ell
_{4}v,y_{6}\right) _{H_{\rho }^{1}(\Omega )}=0.  \label{e5.17}
\end{eqnarray}%
implies that $y_{1}=y_{2}=y_{3}=y_{4}=y_{5}=y_{6}=0$ almost everywhere in
the domain $Q,$ where $\left( \left\{ y_{1},y_{2},y_{3}\right\} ,\left\{
y_{4},y_{5},y_{6}\right\} \right) \in R(\mathcal{X})^{\perp }.$

By putting $W\in D_{0}(\mathcal{X})$ in (\ref{e5.17}), we have
\begin{equation}
\left( \mathcal{L}_{1}(u,v),y_{1}\right) _{L^{2}(0,T;L_{\rho }^{2}(\Omega
))}+\left( \mathcal{L}_{2}(u,v),y_{4}\right) _{L^{2}(0,T;L_{\rho
}^{2}(\Omega ))}=0,  \label{e5.18}
\end{equation}%
hence proposition 5.1 implies that: $y_{1}=y_{4}=0$. Thus (\ref{e5.17})
implies
\begin{equation}
\left( \ell _{1}u,y_{2}\right) _{H_{\rho }^{1}(\Omega )}+\left( \ell
_{2}u,y_{3}\right) _{H_{\rho }^{1}(\Omega )}+\left( \ell _{3}v,y_{5}\right)
_{H_{\rho }^{1}(\Omega )}+\left( \ell _{4}v,y_{6}\right) _{H_{\rho
}^{1}(\Omega )}=0,~~~\forall W\in D_{0}(\mathcal{X}),  \label{e5.19}
\end{equation}%
The four sets $\ell _{1}u,\ell _{2}u$, $\ell _{3}v,$ and $\ell _{4}v$ are
independent, and the images of the trace operator $\ell _{1},\ell _{2}$, $%
\ell _{3},$ and $\ell _{4}$ are respectively everywhere dense in the Hilbert
spaces $H_{\rho }^{1}(\Omega ),$ then it follows from (\ref{e5.19}), that $%
y_{2}=y_{3}=y_{5}=y_{6}=0$ almost everywhere in $Q.$

\section{The nonlinear system}

We are now in a position to solve the nonlinear system (\ref{e1.1}). Relying
on the results obtained previously, we apply an iterative process to
establish the existence and uniqueness of the weak solution of the nonlinear
system (\ref{e1.1}). If $(u,v)$ is a solution of system (\ref{e1.1}) and $%
(\psi ,\phi )$ is a solution of the homogeneous system
\begin{equation}
\begin{cases}
^{C}\partial _{0t}^{\beta }\psi -\frac{1}{x}\left( x\psi _{x}\right) _{x}-%
\frac{1}{x}\left( x\psi _{x}\right) _{xt}+z_{1}\phi +\psi _{t}=0, \\
^{C}\partial _{0t}^{\gamma }\phi -\frac{1}{x}\left( x\phi _{x}\right) _{x}-%
\frac{1}{x}\left( x\phi _{x}\right) _{xt}+z_{2}\psi +\phi _{t}=0, \\
\psi (x,0)=\varphi _{1}(x),~~~\psi _{t}(x,0)=\varphi _{2}(x), \\
\phi (x,0)=\psi _{1}(x),~~~\phi _{t}(x,0)=\psi _{2}(x), \\
\psi _{x}(b,t)=0,~~~\phi _{x}(b,t)=0,~~~\int\limits_{0}^{b}x\psi
dx=0,~~~\int\limits_{0}^{b}x\phi dx=0,%
\end{cases}
\label{e6.1}
\end{equation}%
then $(U,V)=(u-\psi ,v-\phi )$ is a solution of the system
\begin{equation}
\left\{
\begin{array}{c}
^{C}\partial _{0t}^{\beta }U-\frac{1}{x}\left( xU_{x}\right) _{x}-\frac{1}{x}%
\left( xU_{x}\right) _{xt}+z_{1}V+U_{t}=F\left( x,t,U,V,U_{x},V_{x}\right) ,
\\
^{C}\partial _{0t}^{\gamma }V-\frac{1}{x}\left( xV_{x}\right) _{x}-\frac{1}{x%
}\left( xV_{x}\right) _{xt}+z_{2}U+V_{t}=G\left( x,t,U,V,U_{x},V_{x}\right) ,
\\
U(x,0)=0,~~~U_{t}(x,0)=0,~~V(x,0)=0,~~~V_{t}(x,0)=0, \\
\int\limits_{0}^{b}xUdx=0,~~~\int\limits_{0}^{b}xVdx=0,\text{ }%
U_{x}(b,t)=0,~~~V_{x}(b,t)=0,%
\end{array}%
\right.  \label{e6.2}
\end{equation}%
where%
\begin{equation*}
F\left( x,t,U,V,U_{x},V_{x}\right) =f\left( x,t,U+\psi ,V+\phi ,U_{x}+\psi
_{x},V_{x}+\phi _{x}\right) ,
\end{equation*}%
and%
\begin{equation*}
G\left( x,t,U,V,U_{x},V_{x}\right) =g\left( x,t,U+\psi ,V+\phi ,U_{x}+\psi
_{x},V_{x}+\phi _{x}\right) .
\end{equation*}%
The functions $F$ and $G$ are Lipschitzian functions%
\begin{eqnarray}
&&F(x,t,u_{1},v_{1},w_{1},d_{1})-F(x,t,u_{2},v_{2},w_{2},d_{2})\rvert
\label{e6.3} \\
&\leq &\delta _{1}(\left\vert u_{1}-u_{2}\right\vert +\left\vert
v_{1}-v_{2}\right\vert +\left\vert w_{1}-w_{2}\right\vert +\left\vert
d_{1}-d_{2}\right\vert ),  \notag
\end{eqnarray}%
\begin{eqnarray}
&&\left\vert
G(x,t,u_{1},v_{1},w_{1},d_{1})-G(x,t,u_{2},v_{2},w_{2},d_{2})\right\vert
\label{e6.4} \\
&\leq &\delta _{2}(\left\vert u_{1}-u_{2}\right\vert +\left\vert
v_{1}-v_{2}\right\vert +\left\vert w_{1}-w_{2}\right\vert +\left\vert
d_{1}-d_{2}\right\vert ),  \notag
\end{eqnarray}%
for all $(x,t)\in Q$. \newline
According to Theorem 5.1, system (\ref{e6.1}) has a unique solution that
depends continuously on $(\varphi _{1},\varphi _{2},\psi _{1},\psi _{2})\in
\left( H_{\rho }^{1}(\Omega )\right) ^{4}$.

We must prove that the system (\ref{e6.2}) admits a unique solution.

Suppose that $w,U,V\in C^{2}(Q)$, such that
\begin{equation}
w(x,T)=0,~~~w_{t}(x,T)=0,~~~\int\limits_{0}^{b}xw(x,t)dx=0.  \label{e6.5}
\end{equation}%
Consider the identity%
\begin{eqnarray}
&&\left( \mathcal{L}_{1}(U,V),\mathcal{J}_{x}(\xi w)\right)
_{L^{2}(0,T;L_{\rho }^{2}(\Omega ))}+\left( \mathcal{L}_{2}(U,V),\mathcal{J}%
_{x}(\xi w)\right) _{L^{2}(0,T;L_{\rho }^{2}(\Omega ))}  \notag \\
&=&\left( F,\mathcal{J}_{x}(\xi w)\right) _{L^{2}(0,T;L_{\rho }^{2}(\Omega
))}+\left( G,\mathcal{J}_{x}(\xi w)\right) _{L^{2}(0,T;L_{\rho }^{2}(\Omega
))}.  \label{e6.6*}
\end{eqnarray}%
In light of the above assumptions, we obtain%
\begin{equation}
\left( ^{C}\partial _{0t}^{\beta }U,\mathcal{J}_{x}(\xi w)\right)
_{L^{2}(0,T;L_{\rho }^{2}(\Omega ))}=\left( U,\partial _{tT}^{\beta }\left(
\mathcal{J}_{x}(\xi w)\right) \right) _{L^{2}(0,T;L_{\rho }^{2}(\Omega ))},
\label{e6.7}
\end{equation}%
\begin{equation}
-\left( \dfrac{1}{x}\left( xU_{x}\right) _{x},\mathcal{J}_{x}(\xi w)\right)
_{L^{2}(0,T;L_{\rho }^{2}(\Omega ))}=\left( U_{x},xw\right)
_{L^{2}(0,T;L_{\rho }^{2}(\Omega ))},  \label{e6.8}
\end{equation}%
\begin{equation}
-\left( \dfrac{1}{x}\left( xU_{x}\right) _{xt},\mathcal{J}_{x}(\xi w)\right)
_{L^{2}(0,T;L_{\rho }^{2}(\Omega ))}=-\left( U_{x},xw_{t}\right)
_{L^{2}(0,T;L_{\rho }^{2}(\Omega ))},  \label{e6.9}
\end{equation}%
\begin{equation}
\left( z_{1}V,\mathcal{J}_{x}(\xi w)\right) _{L^{2}(0,T;L_{\rho }^{2}(\Omega
))}=-z_{1}\left( \mathcal{J}_{x}(\xi V),w\right) _{L^{2}(0,T;L_{\rho
}^{2}(\Omega ))},  \label{e6.10}
\end{equation}%
\begin{equation}
\left( U_{t},\mathcal{J}_{x}(\xi w)\right) _{L^{2}(0,T;L_{\rho }^{2}(\Omega
))}=-\left( U,\mathcal{J}_{x}(\xi w_{t})\right) _{L^{2}(0,T;L_{\rho
}^{2}(\Omega ))}  \label{e6.10*}
\end{equation}%
\begin{equation}
\left( F,\mathcal{J}_{x}(\xi w)\right) _{L^{2}(0,T;L_{\rho }^{2}(\Omega
))}=-\left( \mathcal{J}_{x}(\xi F),w\right) _{L^{2}(0,T;L_{\rho }^{2}(\Omega
))}.  \label{e6.11}
\end{equation}%
Using the symmetry in the system, and inserting equations (\ref{e6.7})-(\ref%
{e6.11}) into (\ref{e6.6*}), yields
\begin{eqnarray}
&&\left( U,\partial _{tT}^{\beta }\left( \mathcal{J}_{x}(\xi w)\right)
\right) _{L^{2}(0,T;L_{\rho }^{2}(\Omega ))}+\left( V,\partial _{tT}^{\gamma
}\left( \mathcal{J}_{x}(\xi w)\right) \right) _{L^{2}(0,T;L_{\rho
}^{2}(\Omega ))}+\left( U_{x},xw\right) _{L^{2}(0,T;L_{\rho }^{2}(\Omega ))}
\notag \\
&&+\left( V_{x},xw\right) _{L^{2}(0,T;L_{\rho }^{2}(\Omega ))}-\left(
U_{x},xw_{t}\right) _{L^{2}(0,T;L_{\rho }^{2}(\Omega ))}-\left(
V_{x},xw_{t}\right) _{L^{2}(0,T;L_{\rho }^{2}(\Omega ))}  \notag \\
&&-z_{1}\left( \mathcal{J}_{x}(\xi V),w\right) _{L^{2}(0,T;L_{\rho
}^{2}(\Omega ))}-z_{2}\left( \mathcal{J}_{x}(\xi U),w\right)
_{L^{2}(0,T;L_{\rho }^{2}(\Omega ))}-\left( U,\mathcal{J}_{x}(\xi
w_{t})\right) _{L^{2}(0,T;L_{\rho }^{2}(\Omega ))}  \notag \\
&&-\left( V,\mathcal{J}_{x}(\xi w_{t})\right) _{L^{2}(0,T;L_{\rho
}^{2}(\Omega ))}  \notag \\
&=&\left( F,\mathcal{J}_{x}(\xi w)\right) _{L^{2}(0,T;L_{\rho }^{2}(\Omega
))}+\left( G,\mathcal{J}_{x}(\xi w)\right) _{L^{2}(0,T;L_{\rho }^{2}(\Omega
))}.  \label{e6.12}
\end{eqnarray}%
We write (\ref{e6.12}) in the form
\begin{equation}
A\left( w,U,V\right) =\left( w,\mathcal{J}_{x}(\xi F)\right)
_{L^{2}(0,T;L_{\rho }^{2}(\Omega ))}+\left( w,\mathcal{J}_{x}(\xi G)\right)
_{L^{2}(0,T;L_{\rho }^{2}(\Omega ))},  \label{e6.13}
\end{equation}%
where $A\left( w,U,V\right) $ denotes the left-hand side of (\ref{e6.12}).

\begin{defi}
\label{def6} A function $(U,V)\in (L^{2}(0,T;H_{\ \rho }^{1}(\Omega )))^{2}$
is called a weak solution of problem (\ref{e6.2}) if (\ref{e6.13}) and
conditions $U_{x}(b,t)=0,V_{x}(b,t)=0$ hold.
\end{defi}

We now consider the iterated system
\begin{equation}
\begin{cases}
^{C}\partial _{0t}^{\beta }U^{(n)}+-\frac{1}{x}\left( xU_{x}^{(n)}\right)
_{x}-\frac{1}{x}\left( xU_{x}^{(n)}\right)
_{xt}+z_{1}V^{(n)}+U_{t}^{(n)}=F\left(
x,t,U^{(n-1)},V^{(n-1)},U_{x}^{(n-1)},V_{x}^{(n-1)}\right) , \\
^{C}\partial _{0t}^{\gamma }V^{(n)}+-\frac{1}{x}\left( xV_{x}^{(n)}\right)
_{x}-\frac{1}{x}\left( xV_{x}^{(n)}\right)
_{xt}+z_{2}U^{(n)}+V_{t}^{(n)}=G\left(
x,t,U^{(n-1)},V^{(n-1)},U_{x}^{(n-1)},V_{x}^{(n-1)}\right) , \\
U^{(n)}(x,0)=0,~~~U_{t}^{(n)}(x,0)=0,~~V^{(n)}(x,0)=0,~~~V_{t}^{(n)}(x,0)=0,
\\
\int\limits_{0}^{b}xU^{(n)}dx=0,~~~\int%
\limits_{0}^{b}xV^{(n)}dx=0,~~~U_{x}^{(n)}(b,t)=0,~~~V_{x}^{(n)}(b,t)=0,%
\end{cases}
\label{e6.14}
\end{equation}%
where the iterated sequence $\left\{ U^{(n)},V^{(n)}\right\} _{n\geq 0}$ is
constructed in the following way: Given: $\left( U^{(0)},V^{(0)}\right)
=(0,0)$ and the element $\left( U^{(n-1)},V^{(n-1)}\right) ,$ then for $%
n=1,2,\ldots ,$ we solve the problem (\ref{e6.14}). Accordig to Theorem 5.1,
for fixed $n,$ each problem (\ref{e6.14}) has a unique solution $\left(
U^{(n)},V^{(n)}\right) .$

If we set $\left( \mathcal{U}^{(n)}(x,t),\mathcal{V}^{(n)}(x,t)\right)
=\left( U^{(n+1)}(x,t)-\right. $ $\left.
U^{(n)}(x,t),V^{(n+1)}(x,t)-V^{(n)}(x,t)\right) ,$ then we have the new
problem
\begin{equation}
\begin{cases}
& ^{C}\partial _{0t}^{\beta }\mathcal{U}^{(n)}+-\frac{1}{x}\left( x\mathcal{U%
}_{x}^{(n)}\right) _{x}-\frac{1}{x}\left( x\mathcal{U}_{x}^{(n)}\right)
_{xt}+z_{1}\mathcal{V}^{(n)}+\mathcal{U}_{t}^{(n)}=H_{1}^{(n-1)}\left(
x,t\right) , \\
& ^{C}\partial _{0t}^{\gamma }\mathcal{V}^{(n)}+-\frac{1}{x}\left( x\mathcal{%
V}_{x}^{(n)}\right) _{x}-\frac{1}{x}\left( x\mathcal{V}_{x}^{(n)}\right)
_{xt}+z_{2}\mathcal{U}^{(n)}+\mathcal{V}_{t}^{(n)}=H_{2}^{(n-1)}\left(
x,t\right) , \\
& \mathcal{U}^{(n)}(x,0)=0,~~~\mathcal{U}_{t}^{(n)}(x,0)=0,~~\mathcal{V}%
^{(n)}(x,0)=0,~~~\mathcal{V}_{t}^{(n)}(x,0)=0, \\
& \int\limits_{0}^{b}x\mathcal{U}^{(n)}dx=0,~~~\int\limits_{0}^{b}x\mathcal{V%
}^{(n)}dx=0.~~~\mathcal{U}_{x}^{(n)}(b,t)=0,~~~\mathcal{V}_{x}^{(n)}(b,t)=0,%
\end{cases}
\label{e6.15}
\end{equation}%
where%
\begin{equation}
H_{1}^{(n-1)}\left( x,t\right) =F\left(
x,t,U^{(n)},U_{x}^{(n)},V^{(n)},V_{x}^{(n)}\right) -F\left(
x,t,U^{(n-1)},U_{x}^{(n-1)},V^{(n-1)},V_{x}^{(n-1)}\right) ,  \label{e6.16}
\end{equation}%
\begin{equation}
H_{2}^{(n-1)}\left( x,t\right) =G\left(
x,t,U^{(n)},U_{x}^{(n)},V^{(n)},V_{x}^{(n)}\right) -G\left(
x,t,U^{(n-1)},U_{x}^{(n-1)},V^{(n-1)},V_{x}^{(n-1)}\right) .  \label{e6.17}
\end{equation}

\begin{lem}
\label{lem6.1} Assume that conditions (\ref{e6.3}) and (\ref{e6.4}) hold,
then for the fractional linearized system (\ref{e6.15}), we have the a
priori estimate
\begin{equation}
\lVert \mathcal{U}^{(n)}\rVert _{L^{2}(0,T;H_{\ \rho }^{1}(\Omega
))}^{2}+\lVert \mathcal{V}^{(n)}\rVert _{L^{2}(0,T;H_{\ \rho }^{1}(\Omega
))}^{2}\leq K^{\ast }\left( \lVert \mathcal{U}^{(n-1)}\rVert
_{L^{2}(0,T;H_{\ \rho }^{1}(\Omega ))}^{2}+\lVert \mathcal{V}^{(n-1)}\rVert
_{L^{2}(0,T;H_{\ \rho }^{1}(\Omega ))}^{2}\right) ,  \label{e6.18}
\end{equation}%
where $K^{\ast }$ is a positive constant given by%
\begin{equation}
K^{\ast }=4\mathcal{Y}^{\ast \ast }e^{T\mathcal{Y}^{\ast \ast }}T\left(
\delta _{1}^{2}+\delta _{2}^{2}\right) .  \label{e6.18*}
\end{equation}
\end{lem}

\textbf{Proof:} The consideration of the inner products in $L_{\rho
}^{2}\left( \Omega \right) $ of the PDEs in (\ref{e6.15}) and the
integro-differential operators
\begin{equation*}
\mathcal{M}_{1}\mathcal{U}^{(n)}=^{C}\partial _{0t}^{\beta }\mathcal{U}%
^{(n)}+\mathcal{U}_{t}^{(n)}-\mathcal{J}_{x}^{2}(\xi \mathcal{U}_{t}^{(n)}),%
\text{ }\mathcal{M}_{2}\mathcal{V}^{(n)}=^{C}\partial _{0t}^{\gamma }%
\mathcal{V}^{(n)}+\mathcal{V}_{t}^{(n)}-\mathcal{J}_{x}^{2}(\xi \mathcal{V}%
_{t}^{(n)}),
\end{equation*}%
respectively, gives the equation
\begin{eqnarray}
&&\left( ^{C}\partial _{0t}^{\beta }\mathcal{U}^{(n)},^{C}\partial
_{0t}^{\beta }\mathcal{U}^{(n)}+\mathcal{U}_{t}^{(n)}-\mathcal{J}%
_{x}^{2}(\xi \mathcal{U}_{t}^{(n)})\right) _{{L_{\ \rho }^{2}\left( \Omega
\right) }}-\left( \dfrac{1}{x}\left( x\mathcal{U}_{x}^{(n)}\right)
_{x},^{C}\partial _{0t}^{\beta }\mathcal{U}^{(n)}+\mathcal{U}_{t}^{(n)}-%
\mathcal{J}_{x}^{2}(\xi \mathcal{U}_{t}^{(n)})\right) _{{L_{\ \rho
}^{2}\left( \Omega \right) }}  \notag \\
&&-\left( \dfrac{1}{x}\left( x\mathcal{U}_{x}^{(n)}\right)
_{xt},^{C}\partial _{0t}^{\beta }\mathcal{U}^{(n)}+\mathcal{U}_{t}^{(n)}-%
\mathcal{J}_{x}^{2}(\xi \mathcal{U}_{t}^{(n)})\right) _{{L_{\ \rho
}^{2}\left( \Omega \right) }}+\left( z_{1}\mathcal{V}^{(n)},^{C}\partial
_{0t}^{\beta }\mathcal{U}^{(n)}+\mathcal{U}_{t}^{(n)}-\mathcal{J}%
_{x}^{2}(\xi \mathcal{U}_{t}^{(n)})\right) _{L_{\ \rho }^{2}\left( \Omega
\right) }  \notag \\
&&+\left( \mathcal{U}_{t}^{(n)},^{C}\partial _{0t}^{\beta }\mathcal{U}^{(n)}+%
\mathcal{U}_{t}^{(n)}-\mathcal{J}_{x}^{2}(\xi \mathcal{U}_{t}^{(n)})\right)
_{L_{\ \rho }^{2}\left( \Omega \right) }  \notag \\
&&\left( ^{C}\partial _{0t}^{\gamma }\mathcal{V}^{(n)},^{C}\partial
_{0t}^{\gamma }\mathcal{V}^{(n)}+\mathcal{V}_{t}^{(n)}-\mathcal{J}%
_{x}^{2}(\xi \mathcal{V}_{t}^{(n)})\right) _{{L_{\ \rho }^{2}\left( \Omega
\right) }}-\left( \dfrac{1}{x}\left( x\mathcal{V}_{x}^{(n)}\right)
_{x},^{C}\partial _{0t}^{\gamma }\mathcal{V}^{(n)}+\mathcal{V}_{t}^{(n)}-%
\mathcal{J}_{x}^{2}(\xi \mathcal{V}_{t}^{(n)})\right) _{{L_{\ \rho
}^{2}\left( \Omega \right) }}  \notag \\
&&-\left( \dfrac{1}{x}\left( x\mathcal{V}_{x}^{(n)}\right)
_{xt},^{C}\partial _{0t}^{\gamma }\mathcal{V}^{(n)}+\mathcal{V}_{t}^{(n)}-%
\mathcal{J}_{x}^{2}(\xi \mathcal{V}_{t}^{(n)})\right) _{{L_{\ \rho
}^{2}\left( \Omega \right) }}+\left( z_{2}\mathcal{U}^{(n)},^{C}\partial
_{0t}^{\gamma }\mathcal{V}^{(n)}+\mathcal{V}_{t}^{(n)}-\mathcal{J}%
_{x}^{2}(\xi \mathcal{V}_{t}^{(n)})\right) _{L_{\ \rho }^{2}\left( \Omega
\right) }  \notag \\
&&+\left( \mathcal{V}_{t}^{(n)},^{C}\partial _{0t}^{\gamma }\mathcal{V}%
^{(n)}+\mathcal{V}_{t}^{(n)}-\mathcal{J}_{x}^{2}(\xi \mathcal{V}%
_{t}^{(n)})\right) _{L_{\ \rho }^{2}\left( \Omega \right) }  \notag \\
&=&\left( H_{1}^{(n-1)},^{C}\partial _{0t}^{\beta }\mathcal{U}^{(n)}+%
\mathcal{U}_{t}^{(n)}-\mathcal{J}_{x}^{2}(\xi \mathcal{U}_{t}^{(n)})\right)
_{L_{\rho }^{2}\left( \Omega \right) }+\left( H_{2}^{(n-1)},^{C}\partial
_{0t}^{\gamma }\mathcal{V}^{(n)}+\mathcal{V}_{t}^{(n)}-\mathcal{J}%
_{x}^{2}(\xi \mathcal{V}_{t}^{(n)})\right) _{L_{\rho }^{2}\left( \Omega
\right) }.  \label{e6.19}
\end{eqnarray}%
As in the proof of Theorem 4.1, we obtain%
\begin{eqnarray}
&&\lVert \mathcal{U}^{(n)}\rVert _{\mathcal{W}^{\beta }(Q_{t})}^{2}+\lVert
\mathcal{V}^{(n)}\rVert _{\mathcal{W}^{\gamma }(Q_{t})}^{2}+\lVert \mathcal{U%
}^{(n)}\rVert _{H_{\ \rho }^{1}(\Omega )}^{2}+\lVert \mathcal{V}^{(n)}\rVert
_{H_{\ \rho }^{1}(\Omega )}^{2}  \notag \\
&\leq &\mathcal{Y}^{\ast \ast }e^{T\mathcal{Y}^{\ast \ast }}\left(
\int\limits_{0}^{T}\lVert H_{1}^{(n-1)}\rVert _{L_{\ \rho }^{2}(\Omega
)}^{2}d\tau +\int\limits_{0}^{T}\lVert H_{2}^{(n-1)}\rVert _{L_{\ \rho
}^{2}(\Omega )}^{2}d\tau \right)  \label{e6.20}
\end{eqnarray}%
By dropping the first two terms on the left hand side of (\ref{e6.20}), to
get
\begin{equation}
\lVert \mathcal{U}^{(n)}\rVert _{H_{\ \rho }^{1}(\Omega )}^{2}+\lVert
\mathcal{V}^{(n)}\rVert _{H_{\ \rho }^{1}(\Omega )}^{2}\leq \mathcal{Y}%
^{\ast \ast }e^{T\mathcal{Y}^{\ast \ast }}\left( \int\limits_{0}^{T}\lVert
H_{1}^{(n-1)}\rVert _{L_{\ \rho }^{2}(\Omega )}^{2}d\tau
+\int\limits_{0}^{T}\lVert H_{2}^{(n-1)}\rVert _{L_{\ \rho }^{2}(\Omega
)}^{2}d\tau \right) .  \label{e6.21}
\end{equation}%
According to conditions (\ref{e6.3}) and (\ref{e6.4}), we estimate the
right-hand side of (\ref{e6.21}) to obtain%
\begin{equation}
\int\limits_{0}^{T}\lVert H_{i}^{(n-1)}\rVert _{L_{\ \rho }^{2}(\Omega
)}^{2}d\tau \leq 4\delta _{i}^{2}\left( \lVert \mathcal{U}^{(n-1)}\rVert
_{L^{2}(0,T;H_{\ \rho }^{1}(\Omega ))}^{2}+\lVert \mathcal{V}^{(n-1)}\rVert
_{L^{2}(0,T;H_{\ \rho }^{1}(\Omega ))}^{2}\right) ,~~~i=1,2.  \label{e6.22*}
\end{equation}%
Hence, inequality (\ref{e6.21}) becomes
\begin{equation}
\lVert \mathcal{U}^{(n)}\rVert _{H_{\ \rho }^{1}(\Omega )}^{2}+\lVert
\mathcal{V}^{(n)}\rVert _{H_{\ \rho }^{1}(\Omega )}^{2}\leq 4\mathcal{Y}%
^{\ast \ast }e^{T\mathcal{Y}^{\ast \ast }}\left( \delta _{1}^{2}+\delta
_{2}^{2}\right) \left( \lVert \mathcal{U}^{(n-1)}\rVert _{L^{2}(0,T;H_{\
\rho }^{1}(\Omega ))}^{2}+\lVert \mathcal{V}^{(n-1)}\rVert _{L^{2}(0,T;H_{\
\rho }^{1}(\Omega ))}^{2}\right) ,  \label{e6.23}
\end{equation}%
By integrating both sides of (\ref{e6.23}) with respect to $t$ over the
interval $[0,T]$, we obtain
\begin{equation}
\lVert \mathcal{U}^{(n)}\rVert _{L^{2}(0,T;H_{\ \rho }^{1}(\Omega
))}^{2}+\lVert \mathcal{V}^{(n)}\rVert _{L^{2}(0,T;H_{\ \rho }^{1}(\Omega
))}^{2}\leq K^{\ast }\left( \lVert \mathcal{U}^{(n-1)}\rVert
_{L^{2}(0,T;H_{\ \rho }^{1}(\Omega ))}^{2}+\lVert \mathcal{V}^{(n-1)}\rVert
_{L^{2}(0,T;H_{\ \rho }^{1}(\Omega ))}^{2}\right) .  \label{e6.28}
\end{equation}%
where $K^{\ast }$ is given by (\ref{e6.18*}). This achieves the proof of
Lemma 6.1.

\begin{thm}
\label{thm6.1} Suppose that conditions (\ref{e6.3}), and (\ref{e6.4}) hold,
and $K^{\ast }<1/4$, then the nonlinear fractional system (\ref{e6.2})
admits a weak solution in $L^{2}(0,T;H_{\rho }^{1}(\Omega ))$.
\end{thm}

\textbf{Proof:} From (\ref{e6.28}), we conclude that the series $%
\sum_{n=1}^{\infty }\mathcal{U}^{(n)}$ and $\sum_{n=1}^{\infty }\mathcal{V}%
^{(n)}$ converge if $K^{\ast }<1/4.$

Indeed, inequality (\ref{e6.28}), implies%
\begin{equation}
\lVert \mathcal{U}^{(n)}\rVert _{L^{2}(0,T;H_{\ \rho }^{1}(\Omega ))}\leq
\sqrt{K^{\ast }}\left( \lVert \mathcal{U}^{(n-1)}\rVert _{L^{2}(0,T;H_{\
\rho }^{1}(\Omega ))}^{2}+\lVert \mathcal{V}^{(n-1)}\rVert _{L^{2}(0,T;H_{\
\rho }^{1}(\Omega ))}^{2}\right) ^{1/2},  \label{e6.29}
\end{equation}%
\begin{equation}
\lVert \mathcal{V}^{(n)}\rVert _{L^{2}(0,T;H_{\ \rho }^{1}(\Omega ))}\leq
\sqrt{K^{\ast }}\left( \lVert \mathcal{U}^{(n-1)}\rVert _{L^{2}(0,T;H_{\
\rho }^{1}(\Omega ))}^{2}+\lVert \mathcal{V}^{(n-1)}\rVert _{L^{2}(0,T;H_{\
\rho }^{1}(\Omega ))}^{2}\right) ^{1/2}.  \label{e6.30}
\end{equation}%
It follows from (\ref{e6.29}) and (\ref{e6.30}) that%
\begin{eqnarray}
&&\lVert \mathcal{U}^{(n)}\rVert _{L^{2}(0,T;H_{\ \rho }^{1}(\Omega
))}+\lVert \mathcal{V}^{(n)}\rVert _{L^{2}(0,T;H_{\ \rho }^{1}(\Omega ))}
\notag \\
&\leq &2\sqrt{K^{\ast }}\left( \lVert \mathcal{U}^{(n-1)}\rVert
_{L^{2}(0,T;H_{\ \rho }^{1}(\Omega ))}^{2}+\lVert \mathcal{V}^{(n-1)}\rVert
_{L^{2}(0,T;H_{\ \rho }^{1}(\Omega ))}^{2}\right) ^{1/2}.  \label{e6.31}
\end{eqnarray}%
Now since%
\begin{equation}
\lVert \mathcal{U}^{(n)}+\mathcal{V}^{(n)}\rVert _{L^{2}(0,T;H_{\ \rho
}^{1}(\Omega ))}\leq \lVert \mathcal{U}^{(n)}\rVert _{L^{2}(0,T;H_{\ \rho
}^{1}(\Omega ))}+\lVert \mathcal{V}^{(n)}\rVert _{L^{2}(0,T;H_{\ \rho
}^{1}(\Omega ))},  \label{e6.32}
\end{equation}%
then, we infer from (\ref{e6.31}) and (\ref{e6.32}) that%
\begin{eqnarray}
&&\lVert \mathcal{U}^{(n)}+\mathcal{V}^{(n)}\rVert _{L^{2}(0,T;H_{\ \rho
}^{1}(\Omega ))}  \notag \\
&\leq &2\sqrt{K^{\ast }}\left( \lVert \mathcal{U}^{(n-1)}\rVert
_{L^{2}(0,T;H_{\ \rho }^{1}(\Omega ))}^{2}+\lVert \mathcal{V}^{(n-1)}\rVert
_{L^{2}(0,T;H_{\ \rho }^{1}(\Omega ))}^{2}\right) ^{1/2}  \notag \\
&\leq &2\sqrt{K^{\ast }}\left( \lVert \mathcal{U}^{(n-1)}+\mathcal{V}%
^{(n-1)}\rVert _{L^{2}(0,T;H_{\ \rho }^{1}(\Omega ))}^{2}\right) ^{1/2}
\notag \\
&=&2\sqrt{K^{\ast }}\lVert \mathcal{U}^{(n-1)}+\mathcal{V}^{(n-1)}\rVert
_{L^{2}(0,T;H_{\ \rho }^{1}(\Omega ))}.  \label{e6.33}
\end{eqnarray}%
Inequality (\ref{e6.33}), shows that the series $\sum_{n=1}^{\infty }\left(
\mathcal{U}^{(n)}+\mathcal{V}^{(n)}\right) =$ $\sum_{n=1}^{\infty }\mathcal{U%
}^{(n)}+\sum_{n=1}^{\infty }\mathcal{V}^{(n)}$converges if $K^{\ast }<1/4.$
Since $\left( \mathcal{U}^{(n)},\mathcal{V}^{(n)}\right)
=(U^{(n+1)}-U^{(n)},V^{(n+1)}-V^{(n)})$, then it follows that the sequence $%
\left( {U}^{(n)},{V}^{(n)}\right) _{n\in N}$ \ with ${U}^{(n)},$and ${V}%
^{(n)}$defined by:
\begin{eqnarray}
U^{(n)}(x,t) &=&\sum_{k=0}^{n-1}\mathcal{U}^{(k)}(x,t)+U^{(0)}(x,t)  \notag
\\
&=&\sum_{k=0}^{n-1}\left( {U}^{(k+1)}-{U}^{(k)}\right)
+U^{(0)}(x,t),~~~~n=1,2,.....  \label{e6.34}
\end{eqnarray}%
and
\begin{eqnarray}
V^{(n)}(x,t) &=&\sum_{k=0}^{n-1}\mathcal{V}^{(k)}(x,t)+V^{(0)}(x,t)  \notag
\\
&=&\sum_{k=0}^{n-1}\left( {V}^{(k+1)}-{V}^{(k)}\right)
+V^{(0)}(x,t),~~~~n=1,2,.....  \label{e6.35}
\end{eqnarray}%
converge to an element $\left( {U},{V}\right) \in \left( L^{2}(0,T;H_{\ \rho
}^{1}(\Omega ))\right) ^{2},$ which must be proved that it is a solution of
problem (\ref{e6.2}). In other words, $\left( {U},{V}\right) $ must satisfy (%
\ref{e6.13}), and the Neumann boundary conditions.

From the iterated system (\ref{e6.14}), we have:
\begin{eqnarray}
A\left( w,U^{(n)},V^{(n)}\right) &=&\left( w,\mathcal{J}_{x}\left( \xi
F\left( \xi ,t,U^{(n-1)},U_{\xi }^{(n-1)},V^{(n-1)},V_{\xi }^{(n-1)}\right)
\right) \right) _{L^{2}(0,T;L_{\rho }^{2}(\Omega ))}  \notag \\
&&+\left( w,\mathcal{J}_{x}\left( \xi G\left( \xi ,t,U^{(n-1)},U_{\xi
}^{(n-1)},V^{(n-1)},V_{\xi }^{(n-1)}\right) \right) \right)
_{L^{2}(0,T;L_{\rho }^{2}(\Omega ))}  \label{e6.36}
\end{eqnarray}%
We infer from (\ref{e6.36}) that
\begin{eqnarray}
&&A\left( w,U^{(n)}-U,V^{(n)}-V\right) +A\left( w,U,V\right)  \notag \\
&=&\left( w,\mathcal{J}_{x}\left( \xi F\left( \xi ,t,U^{(n-1)},U_{\xi
}^{(n-1)},V^{(n-1)},V_{\xi }^{(n-1)}\right) \right) -\mathcal{J}_{x}\left(
\xi F\left( \xi ,t,U,U_{\xi },V,V_{\xi }\right) \right) \right)
_{L^{2}(0,T;L_{\rho }^{2}(\Omega ))}  \notag \\
&&+\left( w,\mathcal{J}_{x}\left( \xi G\left( \xi ,t,U^{(n-1)},U_{\xi
}^{(n-1)},V^{(n-1)},V_{\xi }^{(n-1)}\right) \right) -\mathcal{J}_{x}\left(
\xi G\left( \xi ,t,U,U_{\xi },V,V_{\xi }\right) \right) \right)
_{L^{2}(0,T;L_{\rho }^{2}(\Omega ))}  \notag \\
&&+\left( w,\mathcal{J}_{x}\left( \xi F\left( \xi ,t,U,U_{\xi },V,V_{\xi
}\right) \right) \right) _{L^{2}(0,T;L_{\rho }^{2}(\Omega ))}+\left( w,%
\mathcal{J}_{x}\left( \xi G\left( \xi ,t,U,U_{\xi },V,V_{\xi }\right)
\right) \right) _{L^{2}(0,T;L_{\rho }^{2}(\Omega ))}.  \label{e6.37}
\end{eqnarray}%
Now from the FPDEs in (\ref{e6.14}), we obtain
\begin{eqnarray}
&&A\left( w,U^{(n)}-U,V^{(n)}-V\right)  \notag \\
&=&\left( w,^{C}\partial _{0t}^{\beta }\mathcal{J}_{x}\left( \xi \left(
U^{(n)}-U\right) \right) \right) _{L^{2}(0,T;L_{\rho }^{2}(\Omega ))}-\left(
w,\mathcal{J}_{x}\left( \frac{\partial }{\partial \xi }\left( \xi \frac{%
\partial }{\partial \xi }\left( U^{(n)}-U\right) \right) \right) \right)
_{L^{2}(0,T;L_{\rho }^{2}(\Omega ))}  \notag \\
&&-\left( w,\frac{\partial }{\partial t}\mathcal{J}_{x}\left( \frac{\partial
}{\partial \xi }\left( \xi \frac{\partial }{\partial \xi }\left(
U^{(n)}-U\right) \right) \right) \right) _{L^{2}(0,T;L_{\rho }^{2}(\Omega
))}+z_{1}\left( w,\mathcal{J}_{x}\left( \xi \left( V^{(n)}-V\right) \right)
\right) _{L^{2}(0,T;L_{\rho }^{2}(\Omega ))}  \notag \\
&&+\left( w,\frac{\partial }{\partial t}\mathcal{J}_{x}\left( \xi \left(
U^{(n)}-U\right) \right) \right) _{L^{2}(0,T;L_{\rho }^{2}(\Omega ))}+\left(
w,^{C}\partial _{0t}^{\gamma }\mathcal{J}_{x}\left( \xi \left(
V^{(n)}-V\right) \right) \right) _{L^{2}(0,T;L_{\rho }^{2}(\Omega ))}  \notag
\\
&&-\left( w,\mathcal{J}_{x}\left( \frac{\partial }{\partial \xi }\left( \xi
\frac{\partial }{\partial \xi }\left( V^{(n)}-V\right) \right) \right)
\right) _{L^{2}(0,T;L_{\rho }^{2}(\Omega ))}-\left( w,\frac{\partial }{%
\partial t}\mathcal{J}_{x}\left( \frac{\partial }{\partial \xi }\left( \xi
\frac{\partial }{\partial \xi }\left( V^{(n)}-V\right) \right) \right)
\right) _{L^{2}(0,T;L_{\rho }^{2}(\Omega ))}  \notag \\
&&+z_{2}\left( w,\mathcal{J}_{x}\left( \xi \left( U^{(n)}-U\right) \right)
\right) _{L^{2}(0,T;L_{\rho }^{2}(\Omega ))}+\left( w,\frac{\partial }{%
\partial t}\mathcal{J}_{x}\left( \xi \left( V^{(n)}-V\right) \right) \right)
_{L^{2}(0,T;L_{\rho }^{2}(\Omega ))}  \label{e6.38*}
\end{eqnarray}%
Conditions on functions $w,U,V,$ and integration of each term on the
right-hand side of (\ref{e6.38*}), yield
\begin{eqnarray}
&&A\left( w,U^{(n)}-U,V^{(n)}-V\right)  \notag \\
&=&-\left( U^{(n)}-U,~~^{C}\partial _{tT}^{\beta }\mathcal{J}_{x}\left( \xi
w\right) \right) _{L^{2}(0,T;L_{\rho }^{2}(\Omega ))}-\left( \frac{\partial
}{\partial x}\left( U^{(n)}-U\right) ,xw\right) _{L^{2}(0,T;L_{\rho
}^{2}(\Omega ))}  \notag \\
&&+\left( \frac{\partial }{\partial x}\left( U^{(n)}-U\right) ,xw_{t}\right)
_{L^{2}(0,T;L_{\rho }^{2}(\Omega ))}-z_{1}\left( V^{(n)}-V,\mathcal{J}%
_{x}\left( \xi w\right) \right) _{L^{2}(0,T;L_{\rho }^{2}(\Omega ))}  \notag
\\
&&+\left( U^{(n)}-U,\mathcal{J}_{x}\left( \xi w_{t}\right) \right)
_{L^{2}(0,T;L_{\rho }^{2}(\Omega ))}+\left( V^{(n)}-V,\mathcal{J}_{x}\left(
\xi w_{t}\right) \right) _{L^{2}(0,T;L_{\rho }^{2}(\Omega ))}  \notag \\
&&-\left( V^{(n)}-V,~~^{C}\partial _{tT}^{\gamma }\mathcal{J}_{x}\left( \xi
w\right) \right) _{L^{2}(0,T;L_{\rho }^{2}(\Omega ))}-\left( \frac{\partial
}{\partial x}\left( V^{(n)}-V\right) ,xw\right) _{L^{2}(0,T;L_{\rho
}^{2}(\Omega ))}  \notag \\
&&+\left( \frac{\partial }{\partial x}\left( V^{(n)}-V\right) ,xw_{t}\right)
_{L^{2}(0,T;L_{\rho }^{2}(\Omega ))}-z_{2}\left( U^{(n)}-U,\mathcal{J}%
_{x}\left( \xi w\right) \right) _{L^{2}(0,T;L_{\rho }^{2}(\Omega ))}
\label{e6.39}
\end{eqnarray}%
Application of the Cauchy--Schwarz inequality, to the terms on the
right-hand side of (\ref{e6.39}) gives%
\begin{equation}
-\left( U^{(n)}-U,~~^{C}\partial _{tT}^{\beta }\mathcal{J}_{x}\left( \xi
w\right) \right) _{L^{2}(0,T;L_{\rho }^{2}(\Omega ))}\leq \left\Vert
U^{(n)}-U\right\Vert _{L^{2}(0,T;L_{\rho }^{2}(\Omega ))}\left\Vert
^{C}\partial _{tT}^{\beta }\mathcal{J}_{x}\left( \xi w\right) \right\Vert
_{L^{2}(0,T;L_{\rho }^{2}(\Omega ))},  \label{e6.40}
\end{equation}%
\begin{equation}
-\left( \frac{\partial }{\partial x}\left( U^{(n)}-U\right) ,xw\right)
_{L^{2}(0,T;L_{\rho }^{2}(\Omega ))}\leq b\left\Vert \frac{\partial }{%
\partial x}\left( U^{(n)}-U\right) \right\Vert _{L^{2}(0,T;L_{\rho
}^{2}(\Omega ))}\left\Vert w\right\Vert _{L^{2}(0,T;L_{\rho }^{2}(\Omega ))},
\label{e6.41}
\end{equation}%
\begin{equation}
+\left( \frac{\partial }{\partial x}\left( U^{(n)}-U\right) ,xw_{t}\right)
_{L^{2}(0,T;L_{\rho }^{2}(\Omega ))}\leq b\left\Vert \frac{\partial }{%
\partial x}\left( U^{(n)}-U\right) \right\Vert _{L^{2}(0,T;L_{\rho
}^{2}(\Omega ))}\left\Vert w_{t}\right\Vert _{L^{2}(0,T;L_{\rho }^{2}(\Omega
))},  \label{e6.42}
\end{equation}%
\begin{equation}
-z_{1}\left( V^{(n)}-V,\mathcal{J}_{x}\left( \xi w\right) \right)
_{L^{2}(0,T;L_{\rho }^{2}(\Omega ))}\leq z_{1}\left\Vert
V^{(n)}-V\right\Vert _{L^{2}(0,T;L_{\rho }^{2}(\Omega ))}\left\Vert \mathcal{%
J}_{x}\left( \xi w\right) \right\Vert _{L^{2}(0,T;L_{\rho }^{2}(\Omega ))},
\label{e6.43}
\end{equation}%
\begin{equation}
-\left( V^{(n)}-V,~~^{C}\partial _{tT}^{\gamma }\mathcal{J}_{x}\left( \xi
w\right) \right) _{L^{2}(0,T;L^{2}(\Omega ))}\leq \left\Vert
V^{(n)}-V\right\Vert _{L^{2}(0,T;L^{2}(\Omega ))}\left\Vert ^{C}\partial
_{tT}^{\gamma }\mathcal{J}_{x}\left( \xi w\right) \right\Vert
_{L^{2}(0,T;L^{2}(\Omega ))},  \label{e6.44}
\end{equation}%
\begin{equation}
-\left( \frac{\partial }{\partial x}\left( V^{(n)}-V\right) ,xw\right)
_{L^{2}(0,T;L_{\rho }^{2}(\Omega ))}\leq b\left\Vert \frac{\partial }{%
\partial x}\left( V^{(n)}-V\right) \right\Vert _{L^{2}(0,T;L_{\rho
}^{2}(\Omega ))}\left\Vert w\right\Vert _{L^{2}(0,T;L_{\rho }^{2}(\Omega ))},
\label{e6.45}
\end{equation}%
\begin{equation}
+\left( \frac{\partial }{\partial x}\left( V^{(n)}-V\right) ,xw_{t}\right)
_{L^{2}(0,T;L_{\rho }^{2}(\Omega ))}\leq b\left\Vert \frac{\partial }{%
\partial x}\left( V^{(n)}-V\right) \right\Vert _{L^{2}(0,T;L_{\rho
}^{2}(\Omega ))}\left\Vert w_{t}\right\Vert _{L^{2}(0,T;L_{\rho }^{2}(\Omega
))},  \label{e6.46}
\end{equation}%
\begin{equation}
-z_{2}\left( U^{(n)}-U,\mathcal{J}_{x}\left( \xi w\right) \right)
_{L^{2}(0,T;L_{\rho }^{2}(\Omega ))}\leq z_{2}\left\Vert
U^{(n)}-U\right\Vert _{L^{2}(0,T;L_{\rho }^{2}(\Omega ))}\left\Vert \mathcal{%
J}_{x}\left( \xi w\right) \right\Vert _{L^{2}(0,T;L_{\rho }^{2}(\Omega ))},
\label{e6.47}
\end{equation}%
\begin{equation}
\left( U^{(n)}-U,\mathcal{J}_{x}\left( \xi w_{t}\right) \right)
_{L^{2}(0,T;L_{\rho }^{2}(\Omega ))}\leq \left\Vert U^{(n)}-U\right\Vert
_{L^{2}(0,T;L_{\rho }^{2}(\Omega ))}\left\Vert \mathcal{J}_{x}\left( \xi
w_{t}\right) \right\Vert _{L^{2}(0,T;L_{\rho }^{2}(\Omega )),}
\label{e6.47*}
\end{equation}%
\begin{equation}
\left( V^{(n)}-V,\mathcal{J}_{x}\left( \xi w_{t}\right) \right)
_{L^{2}(0,T;L_{\rho }^{2}(\Omega ))}\leq \left\Vert V^{(n)}-V\right\Vert
_{L^{2}(0,T;L_{\rho }^{2}(\Omega ))}\left\Vert \mathcal{J}_{x}\left( \xi
w_{t}\right) \right\Vert _{L^{2}(0,T;L_{\rho }^{2}(\Omega ))}.
\label{e6.47**}
\end{equation}

Combination of equality (\ref{e6.39}) and inequalities (\ref{e6.40})-(\ref%
{e6.47**}), leads to
\begin{eqnarray}
&&A\left( w,U^{(n)}-U,V^{(n)}-V\right)  \notag \\
&\leq &l_{1}\left( \left\Vert U^{(n)}-U\right\Vert _{L^{2}(0,T;H_{\rho
}^{1}(\Omega ))}\right)  \notag \\
&&\times \left(
\begin{array}{c}
\left\Vert \partial _{tT}^{\beta }\mathcal{J}_{x}(\xi w)\right\Vert
_{L^{2}(0,T;L^{2}(\Omega ))}+\left\Vert \mathcal{J}_{x}(\xi w)\right\Vert
_{\ L^{2}(0,T;L^{2}(\Omega ))}+\left\Vert w\right\Vert
_{L^{2}(0,T;L^{2}(\Omega ))} \\
+\left\Vert \mathcal{J}_{x}\left( \xi w_{t}\right) \right\Vert
_{L^{2}(0,T;L_{\rho }^{2}(\Omega )),}+\left\Vert w_{t}\right\Vert
_{L^{2}(0,T;L^{2}(\Omega ))}%
\end{array}%
\right)  \notag \\
&&+l_{2}\left( \left\Vert V^{(n)}-V\right\Vert _{L^{2}(0,T;H_{\rho
}^{1}(\Omega ))}\right)  \label{e6.48} \\
&&\times \left(
\begin{array}{c}
\left\Vert \partial _{tT}^{\gamma }\mathcal{J}_{x}(\xi w)\right\Vert
_{L^{2}(0,T;L^{2}(\Omega ))}+\left\Vert \mathcal{J}_{x}(\xi w)\right\Vert
_{L^{2}(0,T;L^{2}(\Omega ))}+\left\Vert w\right\Vert
_{L^{2}(0,T;L^{2}(\Omega ))} \\
+\left\Vert \mathcal{J}_{x}\left( \xi w_{t}\right) \right\Vert
_{L^{2}(0,T;L_{\rho }^{2}(\Omega )),}+\left\Vert w_{t}\right\Vert
_{L^{2}(0,T;L^{2}(\Omega ))}%
\end{array}%
\right) ,  \notag
\end{eqnarray}%
with
\begin{equation*}
l_{1}=l_{2}=\max \left( 1,b,z_{1},z_{2}\right) .
\end{equation*}%
On the other side, we have%
\begin{eqnarray}
&&\left( w,\mathcal{J}_{x}\left( \xi F\left( \xi ,t,U^{(n-1)},U_{\xi
}^{(n-1)},V^{(n-1)},V_{\xi }^{(n-1)}\right) \right) -\mathcal{J}_{x}\left(
\xi F\left( \xi ,t,U,U_{\xi },V,V_{\xi }\right) \right) \right)
_{L^{2}(0,T;L^{2}(\Omega ))}  \notag \\
&\leq &\frac{\delta _{1}b}{\sqrt{2}}\left\Vert w\right\Vert
_{L^{2}(0,T;L^{2}(\Omega ))}\left( \left\Vert U^{(n)}-U\right\Vert
_{L^{2}(0,T;H_{\rho }^{1}(\Omega ))}+\left\Vert V^{(n)}-V\right\Vert
_{L^{2}(0,T;H_{\rho }^{1}(\Omega ))}\right) ,  \label{e6.49}
\end{eqnarray}%
\begin{eqnarray}
&&\left( w,\mathcal{J}_{x}\left( \xi G\left( \xi ,t,U^{(n-1)},U_{\xi
}^{(n-1)},V^{(n-1)},V_{\xi }^{(n-1)}\right) \right) -\mathcal{J}_{x}\left(
\xi G\left( \xi ,t,U,U_{\xi },V,V_{\xi }\right) \right) \right)
_{L^{2}(0,T;L^{2}(\Omega ))}  \notag \\
&&\frac{\delta _{2}b}{\sqrt{2}}\left\Vert w\right\Vert
_{L^{2}(0,T;L^{2}(\Omega ))}\left( \left\Vert U^{(n)}-U\right\Vert
_{L^{2}(0,T;H_{\rho }^{1}(\Omega ))}+\left\Vert V^{(n)}-V\right\Vert
_{L^{2}(0,T;H_{\rho }^{1}(\Omega ))}\right) .  \label{e6.50}
\end{eqnarray}%
As $n\longrightarrow \infty ,$ it follows from (\ref{e6.48})-(\ref{e6.50}),
and (\ref{e6.37}) that%
\begin{equation*}
A\left( w,U,V\right) =\left( w,\mathcal{J}_{x}(\xi F)\right)
_{L^{2}(0,T;L^{2}(\Omega ))}+\left( w,\mathcal{J}_{x}(\xi G)\right)
_{L^{2}(0,T;L^{2}(\Omega ))}.
\end{equation*}%
To conclude that problem (\ref{e6.2}) admits a weak solution, we must show
that conditions $U_{x}(b,t)=0,~~~V_{x}(b,t)=0$ in (\ref{e6.2}) hold. Since; $%
(U,V)\in \left( L^{2}(0,T;H_{\ \rho }^{1}(\Omega ))\right) ^{2}$, then
\begin{equation*}
\int\limits_{0}^{t}U_{x}(x,s)ds,\int\limits_{0}^{t}V_{x}(x,s)ds\in C(%
\overline{Q}),
\end{equation*}%
from which we conclude that: $U_{x}(b,t)=0$,~~$V_{x}(b,t)=0$ , a.e.

It remains now to prove the uniqueness of solution of system (\ref{e6.2}).

\begin{thm}
\label{thm6.2} If hypotheses (\ref{e6.3}) and (\ref{e6.4}) are satisfied,
then the system (\ref{e6.2}) has only one solution.
\end{thm}

\textbf{Proof:} Suppose that $(U_{1},V_{1}),(U_{2},V_{2})\in \left(
L^{2}(0,T;H_{\ \rho }^{1}(\Omega ))\right) ^{2}$ are two different solutions
of the system (\ref{e6.2}) , then $(\mathcal{U},\mathcal{V}%
)=(U_{1}-U_{2},V_{1}-V_{2})\in \left( L^{2}(0,T;H_{\ \rho }^{1}(\Omega
))\right) ^{2}$, verifies
\begin{equation}
\begin{cases}
^{C}\partial _{0t}^{\beta }\mathcal{U}-\frac{1}{x}\left( x\mathcal{U}%
_{x}\right) _{x}-\frac{1}{x}\left( x\mathcal{U}_{x}\right) _{xt}+\mathcal{V+U%
}_{n}=H_{1}\left( x,t\right) , \\
^{C}\partial _{0t}^{\gamma }\mathcal{V}-\frac{1}{x}\left( x\mathcal{V}%
_{x}\right) _{x}-\frac{1}{x}\left( x\mathcal{V}_{x}\right) _{xt}+\mathcal{U+V%
}_{n}=H_{2}\left( x,t\right) , \\
\mathcal{U}(x,0)=0,~~~\mathcal{U}_{t}(x,0)=0,~~\mathcal{V}(x,0)=0,~~~%
\mathcal{V}_{t}(x,0)=0, \\
\int\limits_{0}^{b}x\mathcal{U}dx=0,~~~\int\limits_{0}^{b}x\mathcal{V}%
dx=0.~~~\mathcal{U}_{x}(b,t)=0,~~~\mathcal{V}_{x}(b,t)=0,%
\end{cases}
\label{e6.51}
\end{equation}%
where%
\begin{equation}
H_{1}\left( x,t\right) =F\left( x,t,U_{1},\left( U_{1}\right)
_{x},V_{1},\left( V_{1}\right) _{x}\right) -F\left( x,t,U_{2},\left(
U_{2}\right) _{x},V_{2},\left( V_{2}\right) _{x}\right) ,  \label{e6.52}
\end{equation}%
\begin{equation}
H_{2}\left( x,t\right) =G\left( x,t,U_{1},\left( U_{1}\right)
_{x},V_{1},\left( V_{1}\right) _{x}\right) -G\left( x,t,U_{2},\left(
U_{2}\right) _{x},V_{2},\left( V_{2}\right) _{x}\right) .  \label{e6.53}
\end{equation}%
We now consider the scalar product in the space $L^{2}(0,T;L^{2}(\Omega ))$
of the PDEs in (\ref{e6.51}) and the differential operators $\mathcal{M}_{1}%
\mathcal{U}=^{C}\partial _{0t}^{\beta }\mathcal{U}+\mathcal{U}_{t}-\mathcal{J%
}_{x}^{2}(\xi \mathcal{U}_{t}),$ $\mathcal{M}_{2}\mathcal{V}=^{C}\partial
_{0t}^{\gamma }\mathcal{V}+\mathcal{V}_{t}-\mathcal{J}_{x}^{2}(\xi \mathcal{V%
}_{t}),$ and follow the same computations as in Lemma 6.1, we obatin
\begin{equation}
\lVert \mathcal{U}\rVert _{L^{2}(0,T;H_{\ \rho }^{1}(\Omega ))}^{2}+\lVert
\mathcal{V}\rVert _{L^{2}(0,T;H_{\ \rho }^{1}(\Omega ))}^{2}\leq K^{\ast
}\left( \lVert \mathcal{U}\rVert _{L^{2}(0,T;H_{\ \rho }^{1}(\Omega
))}^{2}+\lVert \mathcal{V}\rVert _{L^{2}(0,T;H_{\ \rho }^{1}(\Omega
))}^{2}\right) ,  \label{e6.54}
\end{equation}%
where $K^{\ast }$ is the same constant as in Lemma \eqref{lem6.1}. Since $%
K^{\ast }<1/4$, we deduce from (\ref{e6.54}) that
\begin{equation}
(1-K^{\ast })\left( \lVert \mathcal{U}\rVert _{L^{2}(0,T;H_{\ \rho
}^{1}(\Omega ))}^{2}+\lVert \mathcal{V}\rVert _{L^{2}(0,T;H_{\ \rho
}^{1}(\Omega ))}^{2}\right) =0,  \label{e6.55}
\end{equation}%
which implies that $(\mathcal{U},\mathcal{V}%
)=(U_{1}-U_{2},V_{1}-V_{2})=(0,0) $, and hence
\begin{equation*}
U_{1}=U_{2}\in L^{2}(0,T;H_{\ \rho }^{1}(\Omega ))\text{ and }V_{1}=V_{2}\in
L^{2}(0,T;H_{\ \rho }^{1}(\Omega )).
\end{equation*}%
This achieves the proof of Theorem \eqref{thm6.2}.

\textbf{Conclusion}. A Caputo fractional nonlinear pseudohyperbolic system
supplemented by a classical and a nonlocal boundary condition of integral
type is investigated. more precisely, in this research work, we search a
function $u(x,t)$ verifying (1.1). The associated fractional linear problem
is reformulated, and the uniqueness and existence of the strong solutions
are proved in a fractional Sobolev space. A priori bound for the solution is
obtained from which the uniqueness of the solution follows. By using some
density arguments, the solvability of the linear problem is established. To
takle the well posedness of the fractional nonlinear problem, we relied on
the obtained results for the linear fractional system, by applying a certain
iterative process. Our study improves and develops some few existence
results for the fractional initial boundary value problems when using the
method of functional analysis, the so called energy inequality method. We
would like to mention that the application of the used method is a little
complicated while dealing with the posed problem in the presence of the
nonlinear source terms, the fractional terms, the appearence of the
singularity and the nonlocal integral conditions.

\end{document}